\documentclass[fleqn]{article}

\usepackage{arxiv}

\usepackage[utf8]{inputenc}
\usepackage[T1]{fontenc}
\usepackage{hyperref}
\usepackage{url}
\usepackage{booktabs}
\usepackage{amsmath}
\usepackage{amsthm}
\usepackage{amsfonts}
\usepackage{microtype}
\usepackage{mathtools}
\usepackage{enumitem}
\usepackage{graphicx}
\usepackage{minted}
\usepackage[
  backend = biber,
  style = authoryear,
  url = false
]{biblatex}
\usepackage{doi}
\usepackage{caption}
\usepackage{longtable}
\usepackage[scaled = 0.8]{beramono}
\usepackage{lineno}
\usepackage{float}
\usepackage{setspace}

\usepackage{newtx}
\usepackage{bm}

\usepackage{cleveref}
\mathtoolsset{showonlyrefs}

\DeclareMathAlphabet{\mathcal}{OMS}{cmsy}{m}{n}
\SetMathAlphabet{\mathcal}{bold}{OMS}{cmsy}{b}{n}

\newtheorem{theorem}{Theorem}
\newtheorem{lemma}{Lemma}
\newtheorem{corollary}{Corollary}

\theoremstyle{definition}

\numberwithin{definition}{section}
\numberwithin{theorem}{section}
\numberwithin{lemma}{section}
\numberwithin{corollary}{section}

\title{A high-dimensional perspective of assortative mating: definition,
dynamics, and limiting results}

\renewcommand{\shorttitle}{A high-dimensional perspective of assortative mating}

\author{
  \href{https://orcid.org/0009-0008-8689-4901}
  {\includegraphics[scale=0.06]{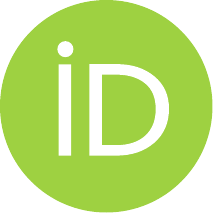}\hspace{1mm}Kári Hlynsson} \\
  Department of Mathematics\\
  Imperial College London\\
  London, SW7 2AZ, UK \\
  \texttt{kari.hlynsson26@imperial.ac.uk}
}

\date{}

\newcommand{\matecor}{R_{\sigma \times \varphi}}
\newcommand{\matecov}{\Sigma_{\sigma \times \varphi}}

\newcommand{\symmatecor}{\overline{R}_{\sigma \times \varphi}}
\newcommand{\symmatecov}{\overline{\Sigma}_{\sigma \times \varphi}}
\newcommand{\asmstrength}{\overline{\rho}_{\sigma \times \varphi}}

\newcommand{\variance}{\operatorname{var}}
\newcommand{\covariance}{\operatorname{cov}}

\newcommand{\expectation}{\operatorname{E}}
\newcommand{\sgn}{\operatorname{sgn}}

\newcommand{\conmult}{\otimes_{\mathfrak{c}}}
\newcommand{\condiv}{\oslash_{\mathfrak{c}}}

\begin{document}

\maketitle
\vspace{-0.75cm}
\begin{abstract}
  Assortative mating, or the systematic (dis)similarity of traits between mate
  pairs, has wide-ranging consequences. Chiefly, it modulates genetic variance
  in offspring due to the coupling of inherited trait-associated loci
  in parents. Over many generations, this buildup becomes appreciable at the
  trait level and may confound statistical genetics techniques, such as
  association analyses and heritability estimation. Theoretical
  characterisations of this effect are well-developed in the univariate case,
  where mates assort on a single trait pair. However, empirical studies
  increasingly show that assortment is not only multivariate, but
  high-dimensional, in that the correlation structure exhibited between mates
  cannot be represented by a single or few linear combinations of the
  constituent traits. In spite of this, the mathematical theory of assortment
  has not scaled with awareness of this fact. In this paper such a
  generalisation is developed. We begin by deriving sufficient conditions under
  which an estimator is invariant to assortment, followed by studying existence
  conditions and limiting behaviour of equilibria under assortment. Finally,
  univariate results involving trait heritability and consanguineal genetic
  correlation are generalised to the multivariate regime.
\end{abstract}


\section{Introduction}
A long-standing and increasingly recognised confounder in the statistical
genetics discipline is the non-random choice of mate in the human population.
This phenomenon, observed as correlation between the traits of spouse pairs,
inflates genetic variance at trait-associated loci in offspring through the
coupled inheritance of these same loci in parents, which becomes appreciable at
the trait level. Compounding through generations, this excess variance and
covariance between loci leads to biased estimates of trait heritability,
genetic correlation within the individual and between relatives, and more.
Importantly, the effects of assortment are not diminished with increased sample
size \parencite{young_solving_2019}, and cannot be fully accounted for through
genetic principal components, which are often used to treat ancestral effects
\parencite{border_assortative_2022}. The theoretical corpus of assortment in the
univariate setting, where mates assort on a single trait pair, is replete and
well formalised, starting from seminal inquiries by
\cite{fisher_xvcorrelation_1919} and \cite{wright_systems_1921}, with subsequent
refinements by \cite{crow_effect_1968}, \cite{crow_introduction_1970},
\cite{feldman_evolution_1977}, \cite{bulmer_mathematical_1980}, and
\cite{nagylaki_assortative_1982}. A partial characterisation of the multivariate
regime is given by \cite{karlin_models_1979} and \cite{karlin_models_1979-1},
although these works are less recognised.

Modern treatments are firmly entrenched in the latter setting, increasingly
augmented by computational methods and simulation \parencite{hartwig_bias_2018,
  border_assortative_2022, border_simple_2024, fan_spousal_2025,
tan_family-gwas_2025}. The need for this change of environment is clear;
assortment has been demonstrated in several surveys to be not only prevalent,
but high-dimensional and the associated spousal cross-correlation to be complex
\parencite{nordsletten_patterns_2016, yengo_imprint_2018, fan_spousal_2025,
border_cross-trait_2022}, ranging from anthropometric, psychiatric, behavioural,
to socioeconomic domains. However, to the best knowledge of the author, this
treatment has not been extended equally to theory. The aim of this paper is to
generalise some of these canonical results to the multivariate regime.

This work is structured as follows.\@ In \Cref{s:methods}, a survey of
notation, machinery, and prior work required to support our efforts is
undertaken. \Cref{ss:results-defining-asm} follows with a definition of
assortment in the high-dimensional regime, alongside characterisation of
\emph{feasibility} and \emph{modality} of assorting regimes; feasibility
describes which assorting structures are feasible, whereas modality assigns them
a positive (assortative), negative (disassortative), random (panmictic), or
indefinite character. The last of these is unprecedented in the univariate
setting. \Cref{ss:results-invariance-to-asm} details a derivation of the
invariance of constant-size, or \emph{atomic}, segments of the genome to
assortment under the infinitesimal polygenic model, which appears as a
convenient utility in subsequent treatment. Following this, the limiting
trajectory and dynamics of the genetic variance matrix are studied in
\Cref{ss:results-dynamics-and-limiting-behaviour}. Under a regularity condition
which is sufficient but not necessary, this recursion tends to a limit which is
used to recover the equilibrium forms of a generalised form of trait
heritability, which under suitable conditions yields an exact recovery of
canonical forms by \cite{fisher_xvcorrelation_1919} and
\cite{gimelfarb_quantitative_1984}. At equilibrium, a derivation of the
covariance between relatives is given in \Cref{ss:covariance-of-relatives}.
Lastly, these identities are validated in \Cref{ss:numerical-study} using
simulation methods.

\section{Methods}
\label{s:methods}

\subsection{Mathematical model and definitions}
\label{ss:mathematical-model-and-definitions}
We study an assorting, outbreeding population numbering \(n \in
\mathbf{N}^\ast\) individuals descended from a randomly mating population of
unrelated founders. Each individual expresses \(p \in \mathbf{N}^\ast\)
quantitative, continuous traits which can be decomposed as the sum of an
additive genetic component and a contribution of the environment. It is further
assumed that the supporting loci for these traits are not sex-linked.

\paragraph{Genotype and phenotype.}
\label{sss:genotype-and-phenotype}
Mathematically, the \(m\) autosomal genetic loci supported in the population are
the collection \(\mathcal G = \{G_{1}, \ldots, G_{m}\}\). Further write
\(\Omega_m = \{0, 1, 2\}^m\) to denote the set of configurations of these loci.
We may also take a subcollection \(\mathcal G^s \subseteq \mathcal G\), in which
case \(\Omega_m^s\) is used to denote configurations in \(\Omega_m\) with
respect to only those loci contained in \(\mathcal G^s\). In particular, such a
collection is said to be \emph{atomic} if its size does not grow with the number
of loci, or mathematically, \(|\mathcal G^s| = O(1)\). Objects in \(\Omega_m\)
and \(\Omega_m^s\) are represented by \(\bm g\) and \(\bm g^s\), or \(\bm G\)
and \(\bm G^s\) if understood as random variables. The notation \(\mu\) is used
to denote the genotype probability distribution, and \(\mu^s\) that of the local
segment \(\mathcal G^s\).

Phenotypes are modelled as the random variables \(X_1, \ldots, X_p\) which we
shall mostly engage with at the vector level \(\bm X = (X_1, \ldots, X_p)\),
taking values in \(\mathbf{R}^p\) and decomposing into an additive genetic
component \(\bm \Gamma\) combined with an environmental component \(\bm E\),
writing \(\bm X = \bm \Gamma + \bm E = \bm G\Lambda + \bm E\), where \(\bm G \in
\Omega_m\), \(\Lambda = (\bm \lambda_1, \ldots, \bm \lambda_p) \in
\mathbf{R}^{m \times p}\) is an effect matrix whose column vector \(\bm\lambda_j
\in \mathbf{R}^m\) designate the causal loci of trait \(j\), and \(\bm E \sim
N(\bm 0, \Sigma_E)\). Indirect components, which involve indirect transmission
of parental genotypes, such as through nurture \parencite{kong_nature_2018} or
cultural transmission \parencite{cavalli-sforza_cultural_1981}, are not
considered here. The matrix \(\Sigma_E\) is the covariance matrix for the
environmental components and is assumed to be non-singular and positive
definite. Effects of trait-associated loci are assumed to be uniform in the
sense that the marginal variance contribution of the locus to the genetic
component is \(m^{-1}\), equivalent to \(\bm\lambda = O(m^{-1/2})\) in the
entrywise sense. The set of loci with non-zero entries in \(\bm\lambda_j\) are
denoted using \(\mathcal G_j\), so \(m_j = |\mathcal G_j|\). It is further
assumed that \(m_j = O(m)\) so that each trait scales with the total number of
loci considered, so that by taking \(m\) to be very large, we are assured that
all traits under consideration are infinitesimally polygenic.

To establish locality results, it will be very convenient to restrict the space
of admissible trait values from Euclidean space \(\mathbf{R}^p\) to some compact
subset \(\mathcal X\). Such an assumption is not as untenable as it may
initially seem; after all, despite the presenting mathematical model, trait
values are usually understood to have some biological floor and ceiling.\@ Two
alternative motivations for this restriction are presented here. On the one
hand, if each constituent trait \(X_j\) can be assigned some compact set of
admissible values \(\mathcal X_j \subset \mathbf{R}\) from a biological
standpoint, we may take the domain as the Cartesian product of these sets,
forming a compact set. If one does not wish to make such an assumption, then an
alternate recourse is defining the domain as the high probability density region
\(\Pr(\bm X \in \mathcal X) \geq 1 - \delta\) for some \(\delta > 0\). This
region may be represented as the preimage of a closed interval in \(\mathbf{R}\)
under a continuous function, and is thus also compact. We shall employ the
latter formalisation. A reader who does not wish to abide by any such assumption
may instead view the asymptotic results presented in this paper as holding with
extremely high probability when \(\delta\) is taken to be small enough and the
number of trait-supporting loci to be sufficiently large or infinite.

\paragraph{Modelling assortment.}\label{ss:modelling-assortment}
Throughout, the symbols \(\sigma\) and \(\varphi\) will be used to denote
spousal or parental traits or genotypes with respect to one or the other sex, in
whichever order. To model assortment, we employ a cross-correlation matrix
\(\matecor \in [-1, 1]^{p \times p}\), the \emph{assorting matrix}. The
cross-covariance is then \(\matecov = \Sigma_{X, \sigma}^{1/2} \matecor
\Sigma_{X, \varphi}^{1/2}\) where \(A^{1/2}\) is used to denote the principal
square root of the positive semi-definite matrix \(A\). Under this model, the
trait vector of either sex can be decomposed as a sum of their expected trait
vector predicated on their mate's traits and independent deviation from this
preference, so that
\begin{equation}\label{eqn:probabilistic-assortment}
  \bm X_{\varphi}
  = \matecov \Sigma_X^{-1} \bm X_{\sigma} + \bm \varepsilon_{\varphi},
\end{equation}
where \(\bm \varepsilon_{\varphi} \sim N(\bm 0_p, \Sigma_{X} -
\matecov^\top \Sigma_X^{-1} \matecov)\). With this model, the desired
cross-covariance structure in mates is recovered while preserving the marginal
variance-covariance structure of both sex classes. The symmetrised form of this
operator, \(\symmatecov = \frac{1}{2} \big(\matecov + \matecov^\top\big)\), is
more amenable to analysis and features extensively. Similarly, \(\symmatecor\)
is used to express this same derivative of the cross-correlation matrix
\(\matecor\). The feasibility and characterisation of the assorting matrix is
discussed in \Cref{ss:results-defining-asm}. Throughout, primary phenotype
assortment, or the independence of spousal genotypes when conditioned on
phenotype, is assumed.

\paragraph{Time and ancestry.}
Time-dependence of model components is indicated with a superscript \((t)\),
where \(t \in \mathbf{N}\) is the index of the generation. Objects such as the
expectation and covariance are often indexed with a superscript to indicate to
which generation the quantity in question corresponds. Importantly, within the
initial panmictic generation \(t = 0\), \(\mu^{(0)}\) is the binomial
distribution that holds in a population at Hardy-Weinberg equilibrium, and
\(\Sigma_G^{(0)}\) is the founder genetic variance-covariance matrix in
\(\mathbf{R}^{m \times m}\).

To represent relatedness of individuals we introduce the directed pedigree graph
\(\mathcal P(t)\), whose vertices are the individuals of generations \(0,
\ldots, t\). Two individuals share an edge if one, say \(j\), is of generation
\(t - 1\) and the other, \(k\), of \(t\) and \(j\) is the parent of \(k\). With
this machinery, relatedness between individuals \(j\) and \(k\) may be
identified with the existence of a path in the undirected version of this graph
connecting the two individuals in question through a most recent common ancestor
(MRCA). This MRCA may be a spouse pair, in which case the individuals are
descendants of siblings, or siblings themselves. Along such a path, write
\(a^d(j)\) with \(1 \leq d \leq d_j\) to denote the individual \(d\) edges
removed from \(j\) on this path, with \(a^{d_j}(j)\) being their MRCA with
\(k\). Of course, this notation corresponds with the iterative application of a
map \(a(j)\) which returns the parent of \(j\) on the path. Similarly, use
\(m^d(j)\) to represent the mate of the individual \(a^d(j)\). Lastly, the sex
of \(a^{d}(j)\) is given by the map \(s(a^{d}(j))\), one of the two tokens
\(\sigma\) or \(\varphi\).

\paragraph{A survey of univariate results.}
A brief overview of the univariate literature is given before proceeding
further. This setting involves a pair of phenotypes \(X = \Gamma + E\) and \(X'
= \Gamma' + E'\), where \(X = X'\) may be the case, where assortment manifests
as the correlation \(r = \mathrm{cor}(X_\sigma, X_\varphi')\). It can be shown
that the shift in additive genetic variance between generations in the trait
\(X\) is
\begin{equation}
  \sigma_{\Gamma, (t + 1)}^2
  = \bigg(1 + \frac{rh_{(t)}^2}{2}\bigg) \sigma_{\Gamma, (t)}^2,
\end{equation}
where \(h_{(t)}^2 = \sigma_{\Gamma, (t)}^2 / \sigma_{X, (t)}^2\) is the
\emph{narrow sense heritability}. In the limit \(t \to \infty\), provided
one exists, the above yields the implicit solution
\begin{equation}
  \frac{\sigma_{\Gamma, (\infty)}^2}{\sigma_{\Gamma, (0)}^2} =
  \frac{1}{1 - rh_{(\infty)}^2}. \label{eqn:wright-limiting-heritability}
\end{equation}
This result was first established by \cite{wright_systems_1921} using a path
tracing model, with subsequent extensions and formalisation by
\cite{crow_effect_1968}, \cite{bulmer_mathematical_1980}, and
\cite{nagylaki_assortative_1982}. \cite{gimelfarb_quantitative_1984} proved
the existence and stability of the limiting heritability, which satisfies the
quadratic equation
\begin{equation}
  r(1 - h_{(0)}^2)\big[h_{(\infty)}^2\big]^2 - h_{(\infty)}^2 + h_{(0)}^2 = 0,
\end{equation}
which has the solution
\begin{equation}
  h_{(\infty)}^2 =
  \begin{dcases}
    \frac{1 - \sqrt{1 - 4rh_{(0)}^2(1 - h_{(0)}^2)}}
    {2r(1 - h_{(0)}^2)}, & r \neq 0, \\
    h_{(0)}^2, & r = 0,
  \end{dcases}
  \label{eqn:gimelfarb_limiting_heritability}
\end{equation}
where the negative branch of the quadratic root is picked so as to ensure
continuity with the panmictic case \(r = 0\) via l'H\^{o}pital's rule. Lastly,
the genetic component covariance of relatives at equilibrium, originating from
\cite{fisher_xvcorrelation_1919}, can be shown to be
\begin{equation}
  \covariance(\Gamma_j, \Gamma_k') = k^{d_j + d_k - 1}
  \sigma_{\Gamma, (\infty)}^2, \label{eqn:kinship-sib-descendants}
\end{equation}
where \(k = (1 + rh^2) / 2\), and \(d_j\) and \(d_k\) denote the number of
steps required to reach the common ancestor of the pair \(j, k\), respectively.
If they are instead related through a single ancestor, such as half-siblings,
this formula is instead
\begin{equation}
  \covariance(\Gamma_j, \Gamma_k') = k^{d_j + d_k} \sigma_{\Gamma, (\infty)}^2.
  \label{eqn:kinship-half-sib-descendants}
\end{equation}

\section{Results}
\label{s:results}

\subsection{Generalising assortative mating to high dimensions}
\label{ss:results-defining-asm}

This subsection investigates the transportation of two attributes in the
univariate setting to the present context: modality, and feasibility. The first
entails whether to speak of assortment in a population as \emph{positive},
\emph{indefinite}, \emph{negative}, or \emph{random}. The second establishes the
class of matrices whose entries lie in \([-1, 1]\) whose elements are admissible
as assorting matrices.

\paragraph{Modes of assortment.}
Since \(\matecor\) need not be symmetric, designating a character of assortment
is not immediate; element-wise classification does not accurately describe the
action induced by this operator, and notions such as definiteness are available
only to symmetric operators. Thus, to bridge between \(\matecor\) and its
symmetrised form, let \(\bm x \in \mathbf{R}^p\), perform the singular value
decomposition \(\matecor = U\Sigma V^\top\), and observe that
\begin{align}
  \bm x^\top \symmatecor \bm x
  = \frac{1}{2} \bm x^\top \big(\matecor + \matecor^\top\big) \bm x
  = \frac{1}{2} \bm x^\top \big(U\Sigma V^\top + V\Sigma U^\top\big) \bm x
  = \sum_{j = 1}^p \sigma_j (\bm x^\top \bm u_j) (\bm x^\top \bm v_j),
  \label{eqn:eigen-vs-singular-cor}
\end{align}
where \(\sigma_j\) is the \(j\)th diagonal entry of \(\Sigma\), and \(\bm u_j\)
and \(\bm v_j\) the \(j\)th column vectors of \(U\) and \(V\), respectively. The
above may be understood as the degree to which the mating regime rewards spouse
pairs whose trait vectors are identical. If \(\bm x\) is a unit vector lying in
the \(k\)th eigenspace \(\mathcal E_k\) of \(\symmatecor\), then the expression
of Equation \eqref{eqn:eigen-vs-singular-cor} reads
\begin{equation}\label{eqn:eigen-vs-sum-of-sing}
  \lambda_k
  = \sum_{j = 1}^p \sigma_j (\bm x^\top \bm u_j) (\bm x^\top \bm v_j).
\end{equation}
Observing that \(\sigma_j \geq 0\), a correspondence between the sign of the
\(k\)th eigenvalue and the sum on the right hand side of Equation
\eqref{eqn:eigen-vs-sum-of-sing} is obtained; if \(\lambda_k < 0\), this implies
that there is at least one summand where \((\bm x^\top \bm u_j) (\bm x^\top \bm
v_j) < 0\), and the contributions of these negative summands must exert a
dominating influence on the sum. In other words, along vectors from \(\mathcal
E_k\), the loading vectors of the sex-respective linear combinations is
"predominantly negative". The opposite interpretation holds for \(\lambda_k >
0\); if \(\lambda_k = 0\), then there is cancellation, meaning that the effects
of assortment averages to zero over canonical trait transformations weighted by
the respective assortment strengths.

If \(\matecor\) is symmetric, then its singular vectors and eigenvectors are the
same up to sign, i.e., \(\bm u_k = \mathrm{sgn}(\lambda_k) \bm v_k\)
and instead of Equation \eqref{eqn:eigen-vs-singular-cor} we get \(\lambda_k =
|\lambda_k| (\bm x^\top \bm u_k)  (\bm x^\top \bm v_k)\). This case admits a
much simpler geometric interpretation. If \(\lambda_k < 0\), then the product
\((\bm x^\top \bm u_k)(\bm x^\top \bm v_k)\) must be negative, and consequently
the vectors geometrically opposed. Conversely, \(\lambda_k > 0\) implies they
must be aligned. In light of the above discussion, it seems natural to define
assortment as \emph{positive} if \(\symmatecor \succ 0\), \emph{negative} if
\(\symmatecor \prec 0\), and \emph{indefinite} if neither of these hold. Of
course, if the assortment matrix is identically zero, then the population is
said to mate randomly.

\paragraph{Feasibility of assorting matrices.}
A description of the set of feasible cross-correlation matrices \(\matecor\) is
now given. To this effect, consider the stacked trait vector \(\bm X_{\sigma +
\varphi} = (\bm X_{\sigma}, \bm X_{\varphi}) \in \mathbf{R}^{2p}\). The variance
of this vector can be expressed as the block matrix
\begin{equation}
  \begin{pmatrix}
    \Sigma_{\sigma} & \matecov \\
    \matecov^\top & \Sigma_{\varphi}
  \end{pmatrix} =
  \begin{pmatrix}
    \Sigma_{X, \sigma} &
    \Sigma_{X, \sigma}^{1/2} \matecor \Sigma_{X, \varphi}^{1/2} \\
    \Sigma_{X, \sigma}^{1/2} \matecor^\top \Sigma_{X, \varphi}^{1/2} &
    \Sigma_{X, \varphi}
  \end{pmatrix},
\end{equation}
which is necessarily a positive semi-definite matrix, hence satisfying the Schur
complement condition
\begin{equation}\label{eqn:asm-feasibility-criterion-svd}
  \Sigma_{X, \varphi}^{1/2} (I - \matecor^\top \matecor)
  \Sigma_{X, \varphi}^{1/2} \succeq 0,
\end{equation}
which requires that \(\matecor^\top \matecor \preceq I\). Equivalently, the
singular values of the assorting matrix must lie in the unit interval
\([0, 1]\). To interpret this constraint, consider the singular value
decomposition of \(\matecor\), which identifies \(p\) orthogonal axes, sometimes
called the \emph{canonical axes}, along which assortment acts independently. The
left and right singular vectors define transformations of the trait space within
each sex class onto these sex-respective axes, and the singular value of each
axis is the (unsigned) correlation between the corresponding pair of whitened,
latent variables. In order for such a decomposition to be valid, the correlation
along these axes must not exceed unity, hence the requirement.

In this section the trait covariance matrices \(\Sigma_{X, \sigma}\) and
\(\Sigma_{X, \varphi}\) were indicated as possibly differing between the sexes;
while such generality will not be considered here, the intention is to
demonstrate that this criterion is irrespective of the covariance structures in
the sexes, and does not require them to be identical.

\subsection{Invariants under assortment}\label{ss:results-invariance-to-asm}
Of central importance to our inquiry is to understand which quantities are
perturbed by assortment. Equally important are those which are not. In this
section an understanding of sufficient conditions for invariance to assortment
are developed. While such a characterisation has been developed in earlier work
\parencite{barton_infinitesimal_2017}, a different one is given here which,
assuming primary phenotypic assortment, establishes conditions for which the
joint probability law of local or atomic, in the sense defined prior, segments
of the genome may be viewed as though derived from a randomly mating population.
Concretely, in \Cref{app:asm-invisible} it is proved that the probability
distribution of an atomic segment \(\mathcal G_s \subset \mathcal G\)
conditional on the trait vector takes the form
\begin{equation}
  \Pr(\bm G^s = \bm g^s \mid \bm X = \bm x) =
  \big(1 + O(m^{-1/2})\big) \Pr(\bm G^s = \bm g^s).
\end{equation}
Thus, in the \(t\)th
generation, where \(t = o(m^{1/2})\), the genotypic probability law
\(\mu_{s}^{(t)}\) satisfies
\begin{equation}
  \mu_{s}^{(t)}(\bm g^s) = \big(1 + O(m^{-1/2})\big) \mu_s^{(0)}(\bm g^s).
\end{equation}
This asymptotic bound of \(t \in \mathbf{N}^\ast\) may be understood as a
restriction of this result to a finite time horizon, which is intuitively the
time period within which accumulated error of iterated approximation of
offspring genotype laws by those of their parents remains negligible. This
restriction on \(t\) may be lifted by taking the limit \(m \to \infty\) first,
making the number of loci infinite, in which case this asymptotic bound is
trivially satisfied. Moreover, in such a setting, this approximation is exact.
In turn, quantities which derive from a marginal or atomic genomic segment are
invariant under assortment in the limiting case. As an example, for a locus
\(G\) with allele frequency \(f\), the variance may be identified as
\begin{equation}
  \variance^{(t)}(G) = \variance^{(0)}(G) = 2f(1 - f),
\end{equation}
seeing as this quantity involves only a single marginal law, which remains
binomially distributed with index \(2\) and event probability \(f\) starting
from Hardy-Weinberg equilibrium. An equally important application of this result
is the segregation variance, a matrix whose entries involve the marginal laws of
two loci, so that
\begin{equation}
  \Sigma_S^{(t)}
  = \Sigma_S^{(0)}
  = \frac{1}{2} \Sigma_G^{(0)}.
\end{equation}
Thus, in what follows, the time dependent notation above is omitted, and we
write \(\Sigma_S\) instead.

\subsection{Dynamics and limiting behaviour}
\label{ss:results-dynamics-and-limiting-behaviour}
Equipped with an understanding of assortment in the multivariate setting, an
investigation into the evolution of an assorting population in time is
undertaken. In concordance with the required protocol to make \(\Sigma_S\) fully
invariant to assortment, these results may be viewed as implicitly taking the
number of loci to be infinite before any limits in the generation index \(t \in
\mathbf{N}^\ast\) are considered. Thus, a more verbose notation would involve
nested limits in \(m\) and \(t\), with the former preceding the latter.

\paragraph{Modulation of genetic variance.}
Consider the offspring genotype vector \(\bm G \in \Omega_m\) and their
parents', \(\bm G_\sigma, \bm G_\varphi \in \Omega_m\) of individuals in
generations \((t + 1)\) and \(t\) generations, respectively. By the law of total
covariance, and using the invariance result of the previous section,
\begin{align}
  \Sigma_G^{(t + 1)}
  &= \expectation[\variance(\bm G) \mid \bm G_\sigma, \bm G_\varphi]
  + \variance(\expectation[\bm G \mid \bm G_\sigma, \bm G_\varphi]) \\
  &= \Sigma_S + \variance(\tfrac{1}{2} \bm G_\sigma + \tfrac{1}{2} \bm G_\varphi
  + \bm \mu) \\
  &= \frac{1}{2} \Sigma_G^{(0)} + \frac{1}{2} \Sigma_G^{(t)}
  + \frac{1}{4} \big(\covariance(\bm G_\sigma, \bm G_\varphi) +
  \covariance(\bm G_\varphi, \bm G_\sigma)\big). \label{eqn:bulmer-first-step}
\end{align}
Write \(\bm G = B_{G \mid X}^{(t)} \bm X + \bm{\varepsilon}\) where \(B_{G
\mid X}^{(t)} = \Sigma_G^{(t)} \Lambda \big[\Sigma_{X}^{(t)}\big]^{-1}\)
predicts \(\bm G\) conditioned on \(\bm X\), and \(\bm\varepsilon \in
\mathbf{R}^m\) is a residual. Then,
\begin{align}
  \mathrm{cov}(\bm{G}_\sigma, \bm{G}_\varphi)
  = \mathrm{cov}(B_{G \mid X}^{(t)} \bm{X}_\sigma + \bm{\varepsilon}_{\sigma},
  B_{G \mid X}^{(t)} \bm{X}_\varphi + \bm{\varepsilon}_{\varphi})
  = B_{G \mid X}^{(t)} \matecov^{(t)} \big[B_{G \mid X}^{(t)}\big]^\top.
\end{align}
By symmetry, \(\covariance(\bm{G}_\varphi, \bm{G}_\sigma) = B_{G \mid X}^{(t)}
\matecov^\top \big[B_{G \mid X}^{(t)}\big]^\top\).
Substituting these into Equation \eqref{eqn:bulmer-first-step} and simplifying,
we have
\begin{equation}\label{eqn:multivariate-bulmer}
  \Sigma_G^{(t + 1)} =
  \frac{1}{2} \Sigma_G^{(0)} + \frac{1}{2} \Sigma_G^{(t)} +
  \frac{1}{2} B_{G \mid X}^{(t)} \symmatecov^{(t)}
  \big[B_{G \mid X}^{(t)}\big]^\top.
\end{equation}
Equation \eqref{eqn:multivariate-bulmer} above is sometimes referred to as
the Bulmer effect in the univariate case and describes how trait correlation
between parents propagates into the genetic variance of offspring.
Unrolling this recursion yields
\begin{align}
  \Sigma_G^{(t + 1)}
  &= \frac{1}{2^{t + 1}} \Sigma_G^{(0)} + \sum_{j = 0}^{t} \frac{1}{2^{j + 1}}
  \Sigma_G^{(0)} + \sum_{j = 0}^{t} \frac{1}{2^{j + 1}} \symmatecov^{(t - j)} \\
  &= \Sigma_G^{(0)} + \sum_{j = 0}^{t} \frac{1}{2^{j + 1}}\symmatecov^{(t - j)}.
\end{align}
Taking the limit \(t \to \infty\), presuming for the moment that a limiting form
exists, gives
\begin{equation}
  \Sigma_G^{(\infty)} = \Sigma_G^{(0)} + \symmatecov^{(\infty)},
\end{equation}
so that the assortment additively inflates or depresses genetic variance
through the symmetrised assorting matrix. Through the view of the above
equation, the preceding discussion centred around the modality of assortment is
made crystal clear; if the assorting matrix is negative, then the limiting
symmetrised assorting covariance above is negative definite, and the result is a
depression in variance. Likewise, if it is positive, an inflation in variance
results from the positivity of the additive term.

\paragraph{Limiting dynamics.}
This section is concerned with the conditions under which the system of Equation
\eqref{eqn:multivariate-bulmer} converges to a limit, although only a
sufficient (but not necessary) characterisation is attempted here. Defining the
map \(\mathcal B\) as the matrix function
\begin{equation}\label{eqn:bulmer-map}
  \mathcal{B}(\Sigma) = \frac{1}{2} \Sigma_G^{(0)} + \frac{1}{2} \Sigma +
  \frac{1}{2} B_{G \mid X}(\Sigma) \symmatecov B_{G \mid X}(\Sigma)^\top,
\end{equation}
the genetic variance matrix sequence can be expressed as \(\Sigma_G^{(t)} =
\mathcal B^{t}(\Sigma_G^{(0)})\), where \(\mathcal B^t\) is used to denote the
application of this map \(t\) times iteratively to its argument. In Appendix
\ref{app:convergence}, it is shown that the matrix iterates can be uniformly
bounded above by a matrix \(\mu I\) where \(\mu > 0\) provided that
the sufficient condition
\begin{equation}
  \asmstrength \theta \|\Sigma_G^{(0)}\| \leq \frac{1}{4},
\end{equation}
holds, where
\(\theta = \lambda_{\max}(\Sigma_\Gamma^{(0)})/\lambda_{\min}(\Sigma_E)\) and
\(\asmstrength = \|\symmatecor\| \leq 1\). Then, using that the map
\(\mathcal{B}(\Sigma)\) is continuous and maps the compact set \(\{0 \preceq
\Sigma \preceq \mu I\}\) into itself, the Brouwer fixed point
theorem yields the existence of a (not necessarily unique) matrix
satisfying \(\mathcal B(\Sigma_G^{(\infty)}) = \Sigma_G^{(\infty)}\). This
condition is likely conservative; in a similar setting,
\cite{karlin_models_1979} showed this limiting form to be unique, with
convergence attained from any initial value, although the model of assortment
considered differs substantially from the present one. Evaluating the form of
Equation \eqref{eqn:bulmer-map} at this fixed point and rearranging,
\begin{equation}\label{eqn:var-fixed-point}
  \Sigma_G^{(\infty)}
  = \Sigma_G^{(0)} + \Sigma_G^{(\infty)} \Lambda \big(\symmatecor \condiv
  \Sigma_X^{(\infty)}\big) \Lambda^\top \Sigma_G^{(\infty)}
\end{equation}
where the notation \(A \condiv B = B^{-1/2} A B^{-1/2}\) and \(A \conmult B =
B^{1/2} A B^{1/2}\) is introduced to simplify exposition. Since the above is
non-linear in the covariance matrix, obtaining a closed form solution is
difficult. However, by careful transformation, the fixed-point equation may be
lifted up to the trait-level, where this endeavour is substantially easier.

\paragraph{Heritability.}
Continuing from Equation \eqref{eqn:var-fixed-point}, multiplying with
\(\Lambda\) from the left and \(\Lambda^\top\) from the right gives the
trait-level identity
\begin{equation}\label{eqn:trait-level-fixed-point}
  \Sigma_\Gamma^{(\infty)} = \Sigma_\Gamma^{(0)} + \Sigma_\Gamma^{(\infty)}
  \left(\symmatecor \condiv \Sigma_X^{(\infty)}\right) \Sigma_\Gamma^{(\infty)}.
\end{equation}
In \Cref{app:convergence} it is shown that the above gives the matrix similarity
relation
\begin{equation}
  \Sigma_\Gamma^{(\infty)} \big[\Sigma_\Gamma^{(0)}\big]^{-1}
  = \symmatecor^{1/2} \big(I - \symmatecor H^{(\infty)}\big)^{-1}
  \symmatecor^{-1/2},
\end{equation}
requiring \(\symmatecor H^{(\infty)} \preceq I\). Assuming pairwise
commutativity of the matrices \((\Sigma_\Gamma^{(0)}, \Sigma_E, \symmatecor)\),
the identity
\begin{equation}
  \frac{\lambda_{\gamma, j}^{(\infty)}}{\lambda_{\gamma, j}^{(0)}}
  =
  \begin{dcases}
    \frac{1}{1 - \overline r_j h_{j, (\infty)}^2} &\overline r_j \neq 0, \\
    1 &\overline r_j = 0,
  \end{dcases}
\end{equation}
may further be identified, which is exactly the form of Equation
\eqref{eqn:wright-limiting-heritability}. This yields the criterion \(r_j h_{j,
(\infty)}^2 < 1\), which may intuitively be understood as the assortment not
being complete along any axis of assortment for which the heredity is absolute.
Lastly, define the \emph{generalised heritability matrix} as \(H^{(t)} =
\Sigma_\Gamma^{(t)} \condiv \Sigma_X^{(t)}\), multiply both sides of
Equation~\eqref{eqn:trait-level-fixed-point} by \(\Sigma_X^{-1/2}\) and
rearrange, yielding
\begin{equation}\label{eqn:heritability-fixed-point}
  H^{(\infty)} \symmatecor H^{(\infty)} - H^{(\infty)} +
  \Theta^{(0)}_{(\infty)} = 0,
\end{equation}
where \(\Theta^{(0)}_{(\infty)} = \Sigma_\Gamma^{(0)} \condiv
\Sigma_X^{(\infty)}\) may be interpreted as the proportion of the heritability
of the traits \(X\) in the limiting population that originates from the
panmictic state. For the solution of this equation, see \Cref{app:convergence}.
Provided that the assorting matrix \(\symmatecor\) is either strictly positive
or negative, or indefinite and one has the pairwise commutativity of
\((\Sigma_\Gamma^{(0)}, \Sigma_E, \symmatecor)\), the solution is
\begin{equation}
  H^{(\infty)} = \frac{1}{2} \sgn(\symmatecor)
  \bigg[\bigg(I - \bigg[I - 4\sgn(\symmatecor)
        \big(\Theta^{(0)}_{(\infty)} \conmult |\symmatecor|\big)
  \bigg]^{1/2}\bigg) \condiv |\symmatecor|\bigg]
  \label{eqn:gimelfarb-noncommuting}
\end{equation}
which requires \(\Theta^{(0)}_{(\infty)} \condiv |\symmatecor| \preceq
\frac{1}{4} I\). In the indefinite and non-commuting case, it is possible that
the discriminant \(\sgn(\symmatecor) (\Theta_{(\infty)}^{(0)} \conmult
|\symmatecor|)\) is not symmetric, and thus a square root involving this
quantity is ill-defined, thus requiring commutativity to be operable. Expanding
on the pairwise commuting case of \((\Sigma_\Gamma^{(0)}, \Sigma_E,
\symmatecor)\), the eigenvalues of the limiting generalised heritability matrix
\(H^{(\infty)}\) can be shown to equal
\begin{equation}
  h_{(\infty), j}^2 =
  \begin{dcases}
    \frac{1 - \sqrt{1 - 4\overline r_j h_{(0), j}^2 (1 - h_{(0), j}^2)}}
    {2\overline r_j(1 - h_{(0), j}^2)}
    &\overline r_j \neq 0, \\
    h_{(0), j}^2 &\overline r_j = 0,
  \end{dcases}\label{eqn:gimelfarb-multivariate}
\end{equation}
which is the form reported in Equation
\eqref{eqn:gimelfarb_limiting_heritability}. Importantly, this solution is valid
for any \(\overline{r}_{j}\) and \(h_{(0), j}^2\).

\subsection{Consanguineal additive genetic covariance}
\label{ss:covariance-of-relatives}
\begin{figure}[t]
  \centering
  \includegraphics[width = 0.95\linewidth]{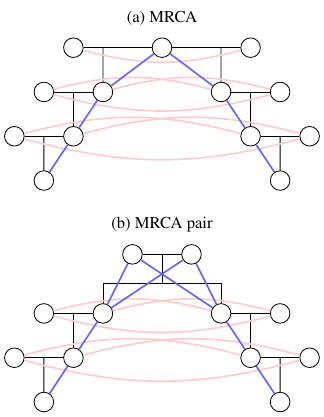}%
  \caption{Schematic diagram of the ancestral lineages of individuals related
    through a single MRCA and a MRCA pair. The blue lines relate the two
    individuals at the bottom of the pedigree through (a) a single most recent
    common ancestor (MRCA), or (b) an MRCA pair. Light-red lines emphasise the
    non-trivial between-spouse covariances in an assorting population along the
  paths connecting relatives.}
  \label{fig:ancestry-diagram}
\end{figure}

The treatment that follows assumes that equilibrium has been reached under
assortment. Suppose that individuals \(j\) and \(k\) of generations \(t_j, t_k
\in \mathbf{N}^\ast\) are related such that there exists a path in the pedigree
\(\mathcal P(t)\) where \(t = \max(t_j, t_k)\) connecting them through either a
single most recent common ancestor (MRCA) or a spouse pair. The former case is
shown in the panel (a) of \Cref{fig:ancestry-diagram} whereas the latter is
shown under (b). Denote by \(d_j\) the length of the end of this path that
connects \(j\) to the MRCA, and similarly \(d_k\) for the opposing part. For
example, in \Cref{fig:ancestry-diagram}, \(d_j = 3\) and \(d_k = 3\) in both
panels. Our inquiry into the additive genetic covariance between \(j\) and \(k\)
starts from the covariance:
\begin{equation}
  \begin{split}
    \mathrm{cov}(\bm\Gamma_j, \bm\Gamma_k)
    &= \expectation
    \!\big[\!\covariance(\bm\Gamma_j, \bm\Gamma_k) \mid \mathcal P(t)\big] +
    \covariance\!\big(\!\expectation[\bm\Gamma_j \mid \mathcal P(t)],
    \expectation[\bm\Gamma_k \mid \mathcal P(t)]\big).
  \end{split}
\end{equation}
Assuming an outbreeding population, the first summand is zero. The expectations
of the covariance term, conditioned on the pedigree \(\mathcal P(t)\), can be
used to write \(\bm\Gamma_j = \frac{1}{2}\bm\Gamma_{a(j)} +
\frac{1}{2}\bm\Gamma_{m(j)} + \bm\mu_j\) and likewise for \(\bm\Gamma_k\). Thus,
\begin{equation}
  \begin{split}
    &\mathrm{cov}(\bm\Gamma_j, \bm\Gamma_k)
    = \frac{1}{4} \bigg(\!\covariance(\bm\Gamma_{a(j)}, \bm\Gamma_{a(k)}) +
      \covariance(\bm\Gamma_{a(j)}, \bm\Gamma_{m(k)}) +
      \covariance(\bm\Gamma_{m(j)}, \bm\Gamma_{a(k)}) +
    \covariance(\bm\Gamma_{m(j)}, \bm\Gamma_{m(k)})\bigg).
  \end{split} \label{eqn:kin-cov-expansion}
\end{equation}
Cross-terms involving individuals \((a(j), m(k))\) and \((m(j), a(k))\), shown
as red lines in \Cref{fig:ancestry-diagram}, cannot be expected to vanish;
indeed, due to mating preferences combined with the genetic similarity of the
ancestors \(a(j)\) and \(a(k)\), we would expect that knowing the trait values
of \(a(j)\) and their relatedness with \(a(k)\) would inform us non-trivially
about their spouse \(m(k)\). Formally, this may be considered by writing
\(\bm\Gamma_{m(j)} = P_{s(j)} \bm\Gamma_{a(j)} + \bm\varepsilon_{m(j)}\), where
\(P_{s(j)}\) is one of the spousal regression matrices
\begin{alignat}{2}
  P_{s(j)} =
  \begin{cases}
    \Sigma_{\Gamma, \sigma \times \varphi}^{(\infty)}
    [\Sigma_\Gamma^{(\infty)}]^{-1} &s(j) = \sigma,
    \\
    [\Sigma_{\Gamma, \sigma \times \varphi}^{(\infty)}]^\top
    [\Sigma_\Gamma^{(\infty)}]^{-1} &s(j) = \varphi.
    \label{eqn:spouse-regression}
  \end{cases}
\end{alignat}
Importantly, in the commuting case we have
\begin{equation}
  P_{s(j)} =
  \begin{cases}
    \symmatecor H^{(\infty)} &s(j) = \sigma, \\
    \symmatecor^\top H^{(\infty)} &s(j) = \varphi.
  \end{cases}
\end{equation}
Using this regression matrix and citing primary phenotypic assortment together
with the outbreeding assumption, the summands of Equation
\eqref{eqn:kin-cov-expansion} can be expressed solely in terms of the ancestral
genetic component covariance, whence
\begin{equation}
  \begin{split}
    \covariance(\bm\Gamma_j, \bm\Gamma_k) =
    &\frac{1}{4} \bigg(
      \covariance(\bm\Gamma_{a(j)}, \bm\Gamma_{a(k)}) +
      \covariance(\bm\Gamma_{a(j)}, \bm\Gamma_{a(k)}) P_{s(k)}^\top + \\
      &\qquad\qquad\qquad P_{s(j)}\covariance(\bm\Gamma_{a(j)},
      \bm\Gamma_{a(k)})
      + P_{s(j)} \covariance(\bm\Gamma_{a(j)}, \bm\Gamma_{a(k)}) P_{s(k)}^\top
    \bigg).
  \end{split}
\end{equation}
This suggests the factorisation
\begin{equation}
  \covariance(\bm\Gamma_j, \bm\Gamma_k) = \left(\frac{I + P_{s(j)}}{2}\right)
  \covariance(\bm\Gamma_{a(j)}, \bm\Gamma_{a(k)})
  \left(\frac{I + P_{s(k)}}{2}\right)^\top.
\end{equation}
\begin{table*}[t]
  \centering
  \begin{tabular*}{\textwidth}{@{\extracolsep{\fill}}lllr@{}}%
    \toprule
    \textbf{\itshape Relationship to self} &
    \textbf{\itshape General case} &
    \textbf{\itshape Commuting, symmetric} &
    \textbf{\itshape Univariate} \\
    \midrule
    Parent &
    \(\mathcal K(j) \Sigma_\Gamma^{(\infty)}\) &
    \(\mathcal K \Sigma_\Gamma^{(\infty)}\) &
    \(k \sigma_{\Gamma, (\infty)}^2\)
    \\
    Grandparent &
    \(\mathcal K^2(j) \Sigma_\Gamma^{(\infty)}\) &
    \(\mathcal K^2 \Sigma_\Gamma^{(\infty)}\) &
    \(k^2 \sigma_{\Gamma, (\infty)}^2\)
    \\
    Half-sibling &
    \(\mathcal K(j) \Sigma_\Gamma^{(\infty)} [\mathcal K(k)]^\top\) &
    \(\mathcal K^2 \Sigma_\Gamma^{(\infty)}\) &
    \(k^2 \sigma_{\Gamma, (\infty)}^2\)
    \\
    Full sibling &
    \(\frac{1}{2} \big[ \Sigma_\Gamma^{(\infty)} +
    \overline{\Sigma}_{\Gamma, \sigma \times \varphi}^{(\infty)} \big]\) &
    \(\mathcal K \Sigma_\Gamma^{(\infty)}\) &
    \(k \sigma_{\Gamma, (\infty)}^2\)
    \\
    First cousin &
    \(\frac{1}{2} \mathcal K(j) \big[ \Sigma_\Gamma^{(\infty)} +
      \overline{\Sigma}_{\Gamma, \sigma \times \varphi}^{(\infty)} \big]
    [\mathcal K(k)]^\top\) &
    \(\mathcal K^3 \Sigma_\Gamma^{(\infty)}\) &
    \(k^3 \sigma_{\Gamma, (\infty)}^2\)
    \\
    Second cousin &
    \(\frac{1}{2} \mathcal K^2(j) \big[ \Sigma_\Gamma^{(\infty)} +
      \overline{\Sigma}_{\Gamma, \sigma \times \varphi}^{(\infty)} \big]
    [\mathcal K^2(k)]^\top\) &
    \(\mathcal K^5 \Sigma_\Gamma^{(\infty)}\) &
    \(k^5 \sigma_{\Gamma, (\infty)}^2\) \\
    Third cousin &
    \(\frac{1}{2} \mathcal K^3(j) \big[ \Sigma_\Gamma^{(\infty)} +
      \overline{\Sigma}_{\Gamma, \sigma \times \varphi}^{(\infty)} \big]
    [\mathcal K^3(k)]^\top\) &
    \(\mathcal K^7 \Sigma_\Gamma^{(\infty)}\) &
    \(k^7 \sigma_{\Gamma, (\infty)}^2\)\\
    \bottomrule
  \end{tabular*}
  \caption{Consanguineal genetic component covariance identities for selected
    familial relationships. The relatives listed in the table are taken as the
  individual \(k\) whereas \(j\) is used as reference.}
  \label{tab:kinship-identities}
\end{table*}
Repeating the above decomposition \(d_j - 1\) and \(d_k - 1\) times on each end
of the lineage, one obtains
\begin{equation}
  \begin{split}
    &\covariance(\bm\Gamma_j, \bm\Gamma_k) =
    \left(\frac{1}{2^{d_j - 1}} \prod_{d = 1}^{d_j - 1} (I + P_{s^d(j)})\right)
    \covariance(\bm\Gamma_{a^{d_j}(j)}, \bm\Gamma_{a^{d_k}(k)})
    \left(\frac{1}{2^{d_k - 1}} \prod_{d = 1}^{d_k - 1}
    (I + P_{s^d(k)})\right)^\top.
  \end{split}
\end{equation}
To simplify, write the matrices on the left and right of the ancestral
covariance as \(\mathcal K^{d_j}(j)\) and \(\mathcal K^{d_k}(k)\), so that, for
example,
\begin{equation}
  \mathcal K^{d_j}(j) = \bigg(\frac{1}{2^{d_j - 1}} \prod_{d = 1}^{d_j - 1}
    \big(
      I + P_{s^d(j)}
  \big)\bigg).
\end{equation}
If \(j\) and \(k\) share a single
MRCA, we may apply the same decomposition once more to obtain the final form
\begin{equation}
  \covariance(\bm\Gamma_j, \bm\Gamma_k)
  = \mathcal{K}^{d_j}(j) \Sigma_\Gamma^{(\infty)}
  \big[\mathcal{K}^{d_k}(k)\big]^\top.
\end{equation}
However, if \(j\) and \(k\) are descendants of (full) siblings or siblings
themselves, then the shared ancestor are their parents, in which case the
boundary term takes the form
\begin{equation}
  \covariance(\bm\Gamma_{a^{d_j}(j)}, \bm\Gamma_{a^{d_k}(k)})
  = \frac{1}{2} (\Sigma_\Gamma^{(\infty)} + \overline{\Sigma}_{\Gamma, \sigma
  \times \varphi}^{(\infty)}).
\end{equation}
Invoking total covariance and citing primary phenotypic assortment to discard
the expected covariance,
\begin{align}
  \Sigma_{\Gamma, \sigma \times \varphi}^{(\infty)}
  &= \expectation[\covariance(
  \bm\Gamma_{\sigma}, \bm\Gamma_{\varphi} \mid \bm X_\sigma, \bm X_\varphi)]
  + \covariance(\expectation[\bm\Gamma_\sigma \mid \bm X_\sigma],
  \expectation[\bm\Gamma_\varphi \mid \bm X_\varphi]) \\
  &= \covariance(\Sigma_\Gamma^{(\infty)} [\Sigma_X^{(\infty)}]^{-1}
    \bm X_\sigma,
  \Sigma_\Gamma^{(\infty)} [\Sigma_X^{(\infty)}]^{-1} \bm X_\varphi)
  \label{eqn:before-simplifying-consang} \\
  &= \Sigma_\Gamma^{(\infty)} (\symmatecor \condiv \Sigma_X^{(\infty)})
  \Sigma_\Gamma^{(\infty)}.
\end{align}
\Cref{tab:kinship-identities} shows the multivariate generalisation compared
with existing univariate results. Despite their somewhat gnarly appearance, this
multivariate form conceals a strong link to the univariate forms, revealed in
the pairwise commuting case considered earlier. There, substituting the
expression of Equation~\eqref{eqn:before-simplifying-consang} into the
covariance allows us to write
\begin{align}
  \covariance(\bm\Gamma_j, \bm\Gamma_k)
  &= \frac{1}{2} \mathcal{K}^{d_j - 1}(j)
  \bigg[\Sigma_\Gamma^{(\infty)} +
  \overline\Sigma_{\Gamma, \sigma \times \varphi}^{(\infty)}\bigg]
  [\mathcal{K}^{d_k - 1}(k)]^\top \\
  &= \frac{1}{2} \mathcal{K}^{d_j - 1}(j)
  \bigg[\Sigma_\Gamma^{(\infty)} +
    \Sigma_{\Gamma}^{(\infty)} [\Sigma_X^{(\infty)}]^{-1}
    \symmatecov
  [\Sigma_X^{(\infty)}]^{-1} \Sigma_{\Gamma}^{(\infty)} \bigg]
  [\mathcal{K}^{d_k - 1}(k)]^\top \\
  &= \frac{1}{2} \mathcal{K}^{d_j - 1}(j)
  \bigg[\Sigma_\Gamma^{(\infty)} +
    \symmatecov
    [H^{(\infty)}]^2
  \bigg]
  [\mathcal{K}^{d_k - 1}(k)]^\top \\
  &= \frac{1}{2} \mathcal{K}^{d_j - 1}(j)
  \bigg[I +
    [\Sigma_\Gamma^{(\infty)}]^{-1}
    \symmatecov
    [H^{(\infty)}]^2
  \bigg] \Sigma_\Gamma^{(\infty)}
  [\mathcal{K}^{d_k - 1}(k)]^\top \\
  &= \frac{1}{2} \mathcal{K}^{d_j - 1}(j)
  \big[I +
    \symmatecor
    H^{(\infty)}
  \big] \Sigma_\Gamma^{(\infty)}
  [\mathcal{K}^{d_k - 1}(k)]^\top \\
  &= \mathcal{K}^{d_j}(j) \Sigma_\Gamma^{(\infty)}
  [\mathcal{K}^{d_k - 1}(k)]^\top,
\end{align}
which closely resembles the univariate form. A further simplification is
afforded if the mating regime is symmetric such that \(\matecor =
\matecor^\top\). In this case \(P_{\sigma} = P_{\varphi}\) and there is no
dependency on sex in the recursion along lineage. In such a setting, the
identities reduce to powers of the matrix \(\mathcal K\) times the boundary
term, and in the commuting case this yields the familiar identity
\begin{equation}
  \mathrm{cov}(\bm\Gamma_j, \bm\Gamma_k) = \mathcal{K}^{d_j + d_k - 1}
  \Sigma_\Gamma^{(\infty)}
\end{equation}
for full siblings or descendants thereof, resembling
Equation~\eqref{eqn:kinship-sib-descendants}, and
\begin{equation}
  \covariance(\bm\Gamma_j, \bm\Gamma_k) = \mathcal{K}^{d_j + d_k}
  \Sigma_\Gamma^{(\infty)}
\end{equation}
for individuals related through a single common ancestor, as in
Equation~\eqref{eqn:kinship-half-sib-descendants}.

\subsection{Numerical study}
\label{ss:numerical-study}

To assess the validity of the theoretical results presented here, analytic
identities were compared against simulated data produced using the forward-time
simulation framework \texttt{knoop} \parencite{hlynsson_knoop-devknoop_2026}.
Results were validated in both the commuting and noncommuting case across all
modalities of assortment, with the assorting matrix symmetric for convenience.
Since the genetic variance identity of Equation~\eqref{eqn:multivariate-bulmer}
involves only the symmetrised form \(\symmatecor\), there is no loss of
generality: indeed, an asymmetric assorting matrix may be exchanged with its
symmetrised counterpart, and the same dynamics will be observed. An exception is
the consanguineal genetic component covariance, which involves the sex-dependent
regression matrices of Equation~\eqref{eqn:spouse-regression} and depends on the
sexes of the individuals appearing on the lineal path connecting relatives. The
initial matrix quantities \(\Sigma_\Gamma^{(0)}\), \(\Sigma_E\), and
\(\matecor\) were generated with respect to a shared orthogonal matrix \(U\), or
mathematically,
\begin{equation}
  \begin{cases}
    \Sigma_\Gamma^{(0)} = UD_\Gamma^{(0)}U^\top, \\
    \Sigma_E = UD_EU^\top, \\
    \matecor = UD_{\sigma \times \varphi} U^\top,
  \end{cases}
\end{equation}
where \(D_{\Gamma}^{(0)}, D_E,\) and \(D_{\sigma \times \varphi}\) are the
spectra of the respective matrices. Non-commuting cases involved generating
distinct orthogonal matrices for each of these respective matrix parameters.
Simulations involved a sex-balanced population of 100,000 individuals, with the
genome consisting of 50,000 unlinked genetic loci in the absence of mutation.
Each simulation scenario was repeated 50 times to measure uncertainty in
estimation. Five traits were simulated, each supported by a total of 10,000 loci
each, so as to provide sufficient polygenicity to approximate the infinitesimal
setting. Locus effect sizes were drawn from a normal distribution and scaled to
the target heritability \(h_j^2\). The complementary environmental component was
drawn from a multivariate normal and scaled such that the resulting traits had
unit variance marginally in the founder generation. In each subsequent
generation, mates were assigned to minimise
\begin{equation}
  \|\tfrac{1}{n}W_\sigma^\top P(s) W_\varphi - \matecor\|_F^2,
\end{equation}
where \(W_\sigma = X_\sigma \Sigma_X^{-1/2}\) and \(W_\varphi = X_\varphi
\Sigma_X^{-1/2}\) are the whitened phenotypes, \(s\) denotes an assignment
of mates, such that \(s(j)\) yields the index of the individual of sex
class \(\varphi\) paired with the individual \(j\) of sex \(\sigma\), and
\(P(s) = (\bm e_{s(1)}, \ldots, \bm e_{s(n)})^\top\) is the permutation matrix
induced by this mate assignment. The matrix \(\matecor\) was drawn from a matrix
normal each generation in accordance with the probabilistic model of
Equation~\eqref{eqn:probabilistic-assortment}, which is important so as not to
omit discrepancies that would be expected in an assorting population of finite
size, despite the true assorting matrix being constant across generations. Once
this assignment was understood in each generation, a number of population
quantities was computed, such as marginal trait heritabilities, and genetic
component covariance within individuals and between relatives, which were
compared against the theory developed here. Across all replicates, a convergence
of genetic component covariance and derivatives was observed, suggesting
uniqueness of the limiting form at the trait level. Moreover, configurations
violating the sufficient criterion \(\asmstrength\theta \|\Sigma^{(0)}\| \leq
1/4\) converged regardless, suggesting that the condition developed here is in
fact conservative.

To contrast simulations, theoretical predictions were computed in R version
4.6.1. In the commuting case, limiting generalised heritability was computed via
Equation~\eqref{eqn:gimelfarb-matrix} in the shared eigenbasis, followed by a
straightforward computation to recover the total and genetic component
covariance matrices of traits at equilibrium, and a recovery of marginal
heritabilities. In the non-commuting case, no closed form is available; instead,
the recursion of Equation~\eqref{eqn:multivariate-bulmer} was iterated to
approximate equilibrium and used to substitute this implicit form in downstream
computation. The equilibrium generalised heritability matrix of
Equation~\eqref{eqn:gimelfarb-noncommuting} was evaluated using this approximate
fixed point.

\begin{figure}[t]
  \centering
  \includegraphics{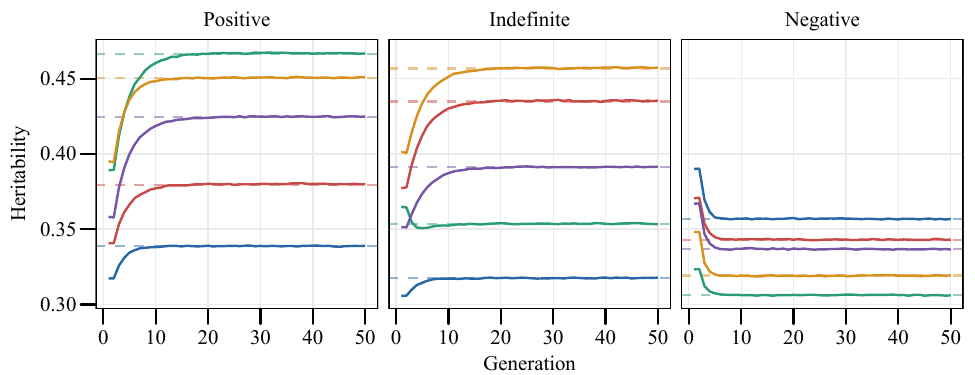}%
  \caption{Evolution of heritability under assortment in the commuting
    case. Simulated heritabilities of five traits under assortment are shown as
    solid lines, with the theoretical limits computed by applying the scalar
    identity of Equation~\eqref{eqn:gimelfarb-multivariate} in the shared
    eigenbasis and recovering the marginal trait heritabilities shown as
  transparent dashed lines, coloured to match the corresponding trait.}
  \label{fig:h2_commuting}
\end{figure}

\begin{figure}[t]
  \centering
  \includegraphics{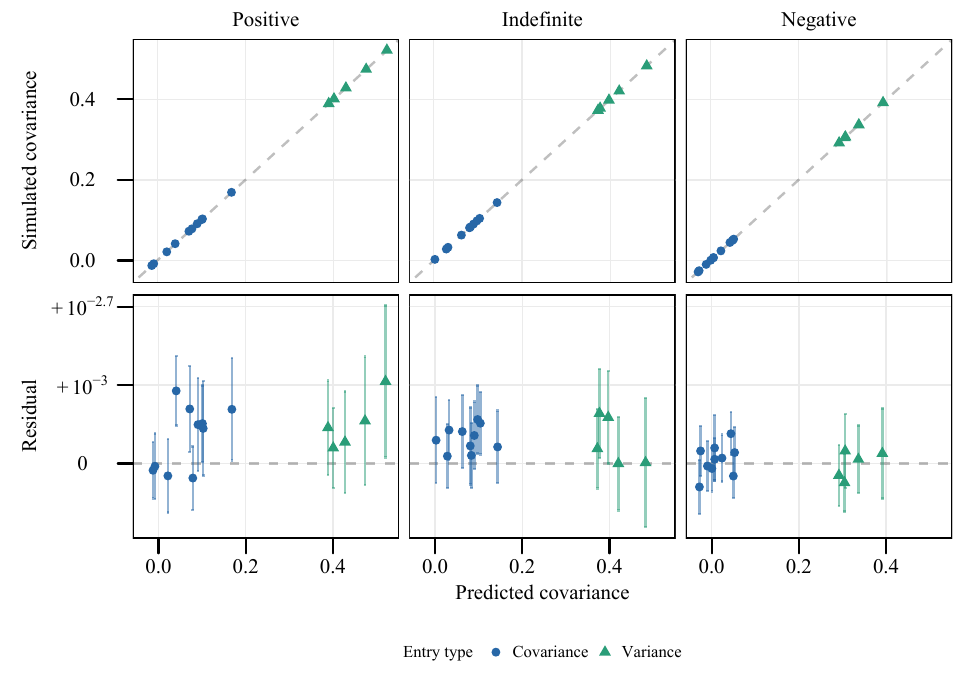}%
  \caption{Equilibrium genetic component covariance under assortment in the
    non-commuting case. The upper panel displays each unique entry of the
    simulated equilibrium covariance matrix against the corresponding entry of
    the predicted matrix obtained by iterating the recursion of
    Equation~\eqref{eqn:multivariate-bulmer} to convergence. The lower panel
    shows the residual between the simulated and predicted entries, with 95\%
  confidence intervals computed across 50 replicates.}
  \label{fig:cov_g_noncommuting}
\end{figure}

\begin{figure}[t]
  \centering
  \includegraphics{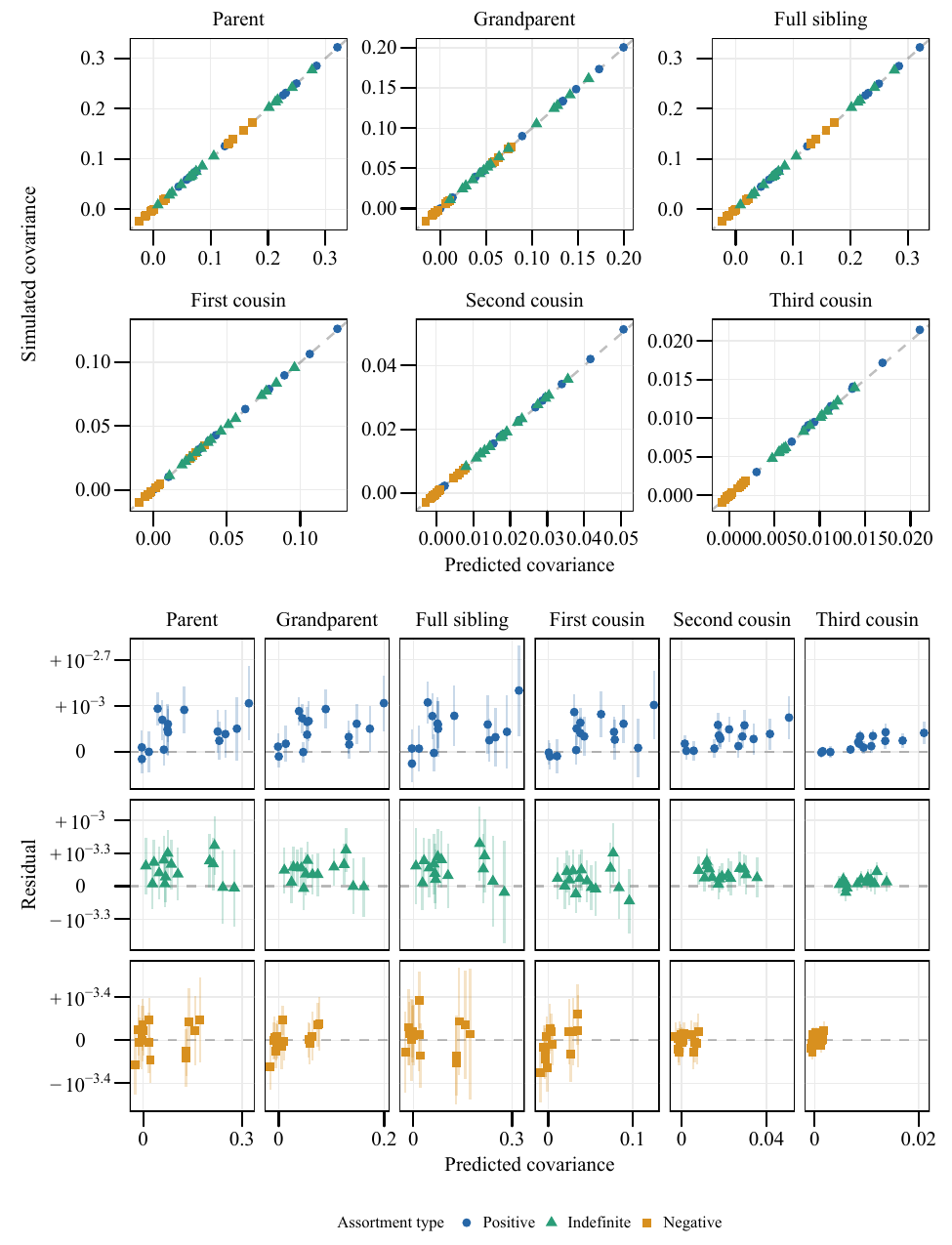}%
  \caption{Consanguineal genetic component covariance under assortment. In the
    upper panel, simulated cross-covariance between relatives are compared
    against the theoretical predictions of \Cref{tab:kinship-identities} at
    equilibrium, for six familial relationships across three assortment
    modalities. Each point represents a unique entry of the corresponding
    covariance matrix. The lower panel displays the residual, along with a 95\%
  confidence interval across 50 replicates.}
  \label{fig:kinship}
\end{figure}

\Cref{fig:h2_commuting,fig:cov_g_noncommuting} shows the evolution of trait
heritabilities and genetic component covariances in simulated data. Convergence
to equilibrium within all traits is observed in approximately 20 generations to
a limit obtained by an evaluation of \Cref{eqn:gimelfarb_limiting_heritability},
followed by algebraic manipulation so as to yield the marginal heritabilities of
the respective traits. \Cref{fig:cov_g_noncommuting} shows both agreement
between simulated and predicted covariance along with entrywise discrepancies in
the noncommuting case, where close agreement is observed. For a majority of
entries, a 95\% confidence interval computed from replicate simulations includes
zero, suggesting this discrepancy is consistent with sampling error, whereas in
some cases this error is elevated away from the origin. This is the case
particularly under positive assortment. However, these errors remain negligible,
of the order of \(10^{-3}\) on traits of initially unit marginal variance.
Lastly, \Cref{fig:kinship} shows the agreement between the predicted and
simulated genetic cross-covariance between relatives and corresponding entrywise
discrepancies. A similar pattern is observed as in
\Cref{fig:cov_g_noncommuting}, with minute errors whose confidence intervals
mostly intersect zero, excepting a handful entries which remain at a slightly
elevated but generally negligible level.

Consistent with prior discussion, the modalities of the assorting matrices in
simulation can be identified with whether the effects of assortment are
inflationary or depressive. In \Cref{fig:h2_commuting}, under a positive
assorting matrix, all trait heritabilities increase monotonically towards
equilibrium, whereas an opposite trend is observed under a negative regime. An
indefinite matrix, by contrast, exhibits both trends; in \Cref{fig:h2_commuting}
a single phenotype experiences depressed variance whereas the rest increase,
although to a lesser extent. This is because the negative latent axes of
assortment load predominantly negatively onto this depressed trait, leading to a
buildup of negative correlation, whereas increasing traits are subject to a
mixed loading of positive and negative components, of which their sum total is
predominantly the former.

Once agreement between theory and simulation was satisfactorily established
through the means described above, a further characterisation of dynamic
properties of the recursion of Equation~\eqref{eqn:multivariate-bulmer} was
sought. To investigate uniqueness, a numerical study was undertaken where a
random configuration model involving components \((\Sigma_G^{(0)}, \Sigma_E,
\matecor)\) was generated, and the fixed point equation for the limiting
genotype covariance implied by Equation~\eqref{eqn:var-fixed-point} solved
numerically via the \texttt{nleqslv} package in R, starting from several random
initial values \(\tilde\Sigma_G^{(0)}\), which, importantly, need not be the
same as \(\Sigma_G^{(0)}\). The components \(\Sigma_G^{(0)}\) and \(\Sigma_E\)
were generated by drawing eigenvalues from a uniform distribution on \([0, 1]\),
followed by separate random orthogonal matrices, assembled as \(A = UD_AU^\top\)
where \(A\) denotes either of the two, \(D_A\) the diagonal matrix of their
respective spectra, and \(U\) the orthogonal matrix for \(A\). The assorting
matrix \(\matecor\) was allowed to be asymmetric by generating two orthogonal
matrices \(U\) and \(V\), such that \(\matecor = U\Sigma V^\top\), where
\(\Sigma\) denotes the desired diagonal matrix of singular values. For
tractability, these investigations were carried out at small scale, with
\(m = 25\); while very small indeed, this is unproblematic since no
infinitesimal
scaling assumptions are required at this level, and the present setup is
sufficient to probe uniqueness of this map.

The above procedure was repeated 1,000 times, with each configuration evaluated
from 50 starting values. Across all configurations, the maximum pairwise
Frobenius distance was of the order \(10^{-7}\), with agreement with the
limiting form obtained through forward iteration of the recursion being
similarly close. Notably, the convergence
criterion \(\asmstrength\theta\|\Sigma_G^{(0)}\| \leq 1/4\) was violated in all
cases, in some cases by several orders of magnitude, yet convergence and
uniqueness were observed throughout. While these results do not rule out the
possibility of a violating state, for which either convergence or uniqueness
fails, they are useful to assess the conservativity of the convergence
criterion discussed prior, and imply that the fixed point of the recursion of
Equation~\eqref{eqn:var-fixed-point} is generally unique.

\section{Discussion}
For some time, assortment has been understood to be a complex and
high-dimensional phenomenon, and has recently become an increasingly focal topic
of research as a force which confounds inference of trait heredity and shared
genetic aetiology. Despite this fact, an extension of canonical univariate
results set forth throughout the 20th century to the multivariate setting has
not been completed, to the best of the author's knowledge. Here such an
endeavour has been undertaken.

It would seem natural to ask whether the univariate equilibrium identities,
such as the limiting heritability and genetic component covariance between
relatives, yield suitable approximations to multivariate regimes. In general,
they do not. Several central quantities used in the analysis here, the
generalised heritability and kinship regression matrix among them, are
indistinguishable from these quantities in the univariate setting, but
substantially different otherwise, which may mislead intuition. This disconnect
is a consequence of dimensionality rather than a particularity of the model
adopted here: when examining the joint system of traits, there are
cross-dependencies which cannot be disregarded in mathematical inquiry. Further,
notions such as the modality of assortment, where a character such as positive,
negative, or indefinite, is attached to an assorting matrix based on its
depressive or inflationary axes along latent trait axes, are ones which may not
be intuitively obvious. For example, one might intuitively attempt to generalise
the definition of assortative and disassortative mating by assigning these
notions entrywise to the assorting matrix, matching the univariate definition.
However, even strong entrywise correlation need not imply an increase or
depression in variance, and thus the former definition is preferable. The aim
of this paper is to develop these notions rigorously so as to supply both
applied and theoretical efforts involving assortment with tools adapted to this
setting.

The simplification afforded under the commutativity of the matrix quantities
\((\Sigma_\Gamma^{(0)}, \Sigma_E, \symmatecor)\) is worthy of remark, both for
its biological significance and for whether such an assumption might be tenable
in practice. In such a setting, traits may be decomposed into orthogonal
eigencomponents with respect to a shared basis, intuitively understood as latent
genetic and environmental traits which combine through linear combinations to
form each of the realised traits. Moreover, the alignment of the assorting
matrix causes this geometry to be preserved intergenerationally, so that the
same latent traits remain foundational to the observed traits, but the
contributions of the latent genetic components are modulated via the spectrum of
the genetic covariance matrix. While appealing from a theoretical perspective,
and perhaps appearing intuitively true under such a framing, this assumption is
not known to hold, and can likely only be verified through empirical inquiry.

Indirect trait effects, such as nurture \parencite{kong_nature_2018} and
cultural transmission effects \parencite{cavalli-sforza_cultural_1981}, have not
been considered here. While the author anticipates such an extension to be
broadly compatible with the strategies employed here, the scope of this
treatment has been set so as to exclude these effects. The reasons for this are
twofold. Firstly, much of the univariate literature which this work extends
incorporates only additive genetic and environmental effects, and thus following
such a convention seems appropriate. Secondly, this restriction of scope has, in
broad strokes, been set so as to match the univariate setting it extends,
providing a foundation from which further generalisation, in both dimensionality
and model complexity, may proceed. These effects may, however, be incorporated
into the framework developed here via the addition of an additive term which
takes the form \((\bm G_\sigma + \bm G_\varphi)\Psi\), where \(\Psi\) is an
effect matrix which designates the contributions of parental genotypes to the
offspring phenotypes, regardless of transmission. If a differential level of
influence between parents is assumed, this addition may be replaced by the
composite term \(\bm G_\sigma \Psi_\sigma + \bm G_\varphi \Psi_\varphi\), where
\(\Psi_\sigma\) and \(\Psi_\varphi\) the indirect transmission effects of each
respective parent. The entries of these matrices may be taken to be proportional
to those of the effect matrix \(\Lambda\), or to be of a different form
altogether. Upon the addition of such a term, a second-order dependency is
injected into the genotype covariance recursion, involving grandparental
genotypes of the offspring. This will be the subject of future work.

Throughout it has been assumed that there are no forces of selection acting on
the population. In particular, fertility has been taken to be independent of the
traits under assortment. It is quite plausible that such dynamics are at play
in reality \parencite{conley_assortative_2016}, however, and these results may
fail to capture these effects. While the reasoning for omitting such effects
remains the same as outlined above, this omission indicates another direction
which would extend this work. Generalisation to sex-linked genotypes would also
meaningfully increase the applicability of the theory presented here to the
applied setting. The assorting matrix has also been assumed constant across
generations, something which may well not be the case in reality
\parencite{schwartz_trends_2013, conley_assortative_2016}. Indeed, it is
plausible that extraneous factors, societal norms being one such example,
modulate the assorting structure in time. Allowing such generality in this
setting would mean replacing the constant assorting matrix with time-indexed
matrices, whereupon the sandwiched term of
Equation~\eqref{eqn:multivariate-bulmer} is replaced with a summation of
conjugacies involving the regression matrices and each temporal assorting
matrix. This form is little different from the one that is presented here in
appearances, but matters of convergence require that a sum over this series
tends to a stable, well-defined limit.

The ancestral composition of the population is a further consideration. In an
admixed population with distinct assorting structure within each subgroup, a
single assorting matrix may be too restrictive. A natural extension would
designate within- and between-group assorting matrices, giving rise to a system
of coupled recursions. Further work might involve such an extension, along with
those detailed above, preceded by empirical inquiry to understand which of these
considerations has meaningful bearing on the model presented here.

\section{Declarations}
\subsection{Acknowledgments}
The author is indebted to Professor Augustine Kong from the University of
Oxford, discussions with whom helped shape the early versions of this work. The
author also extends his gratitude to Professor Heather Battey from Imperial
College London for her support.

\subsection{Code availability}
\texttt{knoop} is distributed under the GNU Public License v3 and accessible
at \url{https://doi.org/10.5281/zenodo.24446743}. Scripts used for the
numerical study are available at \url{https://doi.org/10.5281/zenodo.24447072}.

\printbibliography

\clearpage
\newpage

\appendix

\section{Appendix}

\subsection{Invariants under assortment}\label{app:asm-invisible}

To start, consider the probability law of the trait vector \(\bm X\).
Since \(\bm X = \bm \Gamma + \bm E = \bm G\Lambda + \bm{E}\), the law is given
by the convolution \(f_{\bm X} = f_{\bm\Gamma} \ast f_{\bm E}\). Now, since
\(f_{\bm E}\) is \(N(\bm 0, \Sigma_E)\), the result is a smooth
function, \(f_{\bm X} \in C^\infty(\mathbf{R}^p)\), since the convolution of a
measurable function with a Gaussian is smooth.

\begin{lemma}\label{lemma:state-ratio-asymptotics}
  Let \(\bm x \in \mathcal X\) and \(\bm g^s_1, \bm g^s_2 \in \Omega_m^s\). Then
  \begin{equation}
    \frac{f_{X \mid G^s}(\bm x - \bm g_1^s \Lambda)}{f_{X \mid G^s}(\bm x - \bm
    g_2^s \Lambda)} = 1 + O(m^{-1/2}).
  \end{equation}
\end{lemma}

\begin{proof}
  Write \(\bm\xi = \bm x - \bm g_1^s \Lambda\) and \(\bm \eta = \bm x - \bm
  g_2^s \Lambda\) for convenience. Using that \(f_{X \mid G^s}\) is smooth,
  expand the numerator in a Taylor series around \(\bm \eta\) such that
  \begin{align}
    \frac{f_{X \mid G^s}(\bm\xi)}{f_{X \mid G^s}(\bm\eta)}
    &= 1 + \frac{\nabla f_{X \mid G^s}}{f_{X \mid G^s}}(\bm \eta)
    (\bm g_2^s - \bm g_1^s) \Lambda + O(\|(\bm g_2^s - \bm
    g_1^s)\Lambda\|^2) \\
    &= 1 + \nabla \log f_{X \mid G^s}(\bm\eta)
    (\bm g_2^s - \bm g_1^s) \Lambda + O(\|(\bm g_2^s - \bm
    g_1^s) \Lambda\|^2) \label{eqn:score-statistic}
  \end{align}
  Since \(\bm\eta\) is contained in the compact set \(\mathcal X\) and \(\nabla
  \log f_{X \mid G^s}\) is continuous by virtue of \(f_{X \mid G^s}\) being
  smooth, the former function attains a minimum and maximum value, so we view
  its contribution to Equation \eqref{eqn:score-statistic} as constant. Lastly,
  note that the configurations \(\bm g_1^s\) and \(\bm g_2^s\) differ at most by
  \(|\mathcal G^s| = O(1)\) locus assignments, whereas the effect contribution
  \(\Lambda\) is of the order \(O(m^{-1/2})\). Thus, \(\|(\bm g_2^s - \bm g_1^s)
  \Lambda\| = O(m^{-1/2})\). Revisiting Equation
  \eqref{eqn:score-statistic},
  \begin{align}
    \frac{f_{X \mid G^s}(\bm\xi)}{f_{X \mid G^s}(\bm\eta)}
    &= 1 + \nabla \log f_{X \mid G^s}(\bm\eta) O(m^{-1/2}) +
    O(m^{-1}) \\
    &= 1 + O(m^{-1/2}),
  \end{align}
  as was to be shown.
\end{proof}

It is worth noting that the maximum of the score \(\nabla \log f_{X \mid G^s}\)
over the domain may grow without bound as \(\mathcal X\) expands. This is
unproblematic when the number of loci \(m\) is understood as being sufficiently
large or infinite, but it does indeed require that the number of loci supporting
the traits considered here need to scale with the domain \(\mathcal X\).

\begin{lemma}\label{lemma:probability-asymptotics}
  Let \(\bm g^s \in \Omega_m^s\) and \(\bm x \in \mathcal X\). Then
  \begin{equation}
    \Pr(\bm G^s = \bm g^s \mid \bm X = \bm x)
    = \big(1 + O(m^{-1/2})\big) \Pr(\bm G^s = \bm g^s).
  \end{equation}
\end{lemma}

\begin{proof}
  An application of Bayes' rule gives
  \begin{align}
    \Pr(\bm G^s = \bm g^s \mid \bm X = \bm x)
    &= \frac{\Pr(\bm X = \bm x \mid \bm G^s = \bm g^s)}{\Pr(\bm X = \bm x)}
    \Pr(\bm G^s = \bm g^s) \\
    &= \frac{f_{X \mid G^s}(\bm x - \bm g^s \Lambda)}
    {\sum_{\bm g^s_1 \in \Omega_m^s} f_{G^s}(\bm g_1^s) f_{X \mid G^s}(\bm x -
    \bm g_1^s \Lambda)} \Pr(\bm G^s = \bm g^s) \\
    &=
    \frac{f_{X \mid G^s}(\bm x - \bm g^s \Lambda)}
    {\expectation_{G_s}\!\left[f_{X \mid G^s}(\bm x - \bm G^s \Lambda)\right]}
    \Pr(\bm G^s = \bm g^s) \\
    &= \expectation_{G_s}\!\left[
      \frac{f_{X \mid G^s}(\bm x - \bm G^s \Lambda)}{f_{X \mid G^s}(\bm x -
      \bm g^s \Lambda)}
    \right]^{-1} \Pr(\bm G^s = \bm g^s).
  \end{align}
  Applying Lemma \ref{lemma:state-ratio-asymptotics} and recognising that the
  expectation is a sum ranging over a constant number of elements, the result is
  obtained.
\end{proof}

\begin{theorem}\label{thm:assortment-invariance}
  Let \(\mu^{(t)}_s\) denote the probability law of the genotypes \(\mathcal G^s
  \subset \mathcal G\) of constant size within generation \(t \in \mathbf{N}\),
  where \(\mu^{(0)}_s\) is the law of a randomly mating population. Then, for
  \(t = o(m^{1/2})\) and \(\bm g^s \in \Omega_m^s\),
  \begin{equation}
    \mu_s^{(t)}(\bm g^s) = \big(1 + O(m^{-1/2})\big) \mu^{(0)}_s(\bm
    g^s).
  \end{equation}
\end{theorem}

\begin{proof}
  The proof is by induction, starting from the base case \(t = 1\). By the law
  of total probability,
  \begin{align}
    \mu_s^{(1)}(\bm g^s)
    &= \sum_{\bm g^s_\sigma, \bm g^s_\varphi \in \Omega_m^s}
    \Pr(\bm G^s = \bm g^s \mid \bm G^s_\sigma = \bm g^s_\sigma,
    \bm G^s_\varphi = \bm g^s_\varphi) \Pr(\bm G^s_\sigma = \bm g^s_\sigma,
    \bm G^s_\varphi = \bm g^s_\varphi) \\
    &= \sum_{\bm g^s_\sigma, \bm g^s_\varphi \in \Omega_m^s}
    \psi(\bm g^s \mid \bm g^s_\sigma, \bm g^s_\varphi)
    \Pr(\bm G^s_\sigma = \bm g^s_\sigma, \bm G^s_\varphi = \bm g^s_\varphi)
    \label{eqn:mendelian-geno-law}
  \end{align}
  where \(\psi\) is a Mendelian transmission kernel whose form is independent of
  population structure. Conditioning over the joint law of the parental
  genotypes,
  \begin{align}
    &\Pr(\bm G^s_\sigma = \bm g^s_\sigma, \bm G^s_\varphi = \bm g^s_\varphi) \\
    &\qquad= \sum_{\bm x_\sigma, \bm x_\varphi \in \mathcal X}
    \Pr(\bm G^s_\sigma = \bm g^s_\sigma, \bm G^s_\varphi = \bm g^s_\varphi
    \mid \bm X_\sigma = \bm x_\sigma, \bm X_\varphi = \bm x_\varphi)
    \Pr(\bm X_\sigma = \bm x_\sigma, \bm X_\varphi = \bm x_\varphi) \\
    &\qquad= \sum_{\bm x_\sigma, \bm x_\varphi \in \mathcal X}
    \Pr(\bm G^s_\sigma = \bm g^s_\sigma \mid \bm X_\sigma = \bm x_\sigma)
    \Pr(\bm G^s_\varphi = \bm g^s_\varphi \mid \bm X_\varphi = \bm x_\varphi)
    \Pr(\bm X_\sigma = \bm x_\sigma, \bm X_\varphi = \bm x_\varphi)
    \label{eqn:use-ppa} \\
    &\qquad= \sum_{\bm x_\sigma, \bm x_\varphi}
    \left[\big(1 + O(m^{-1/2})\big) \mu^{(0)}_s(\bm g^s_\sigma)
    \mu^{(0)}_s(\bm g^s_\varphi)\right] \Pr(\bm X_\sigma = \bm x_\sigma,
    \bm X_\varphi = \bm x_\varphi)
    \label{eqn:use-probability-asymptotics} \\
    &\qquad= \big(1 + O(m^{-1/2})\big) \mu^{(0)}_s(\bm g^s_\sigma)
    \mu^{(0)}_s(\bm g^s_\varphi).
  \end{align}%
  Equation \eqref{eqn:use-ppa} invoked primary phenotypic assortment,
  followed by
  an application of Lemma \ref{lemma:probability-asymptotics} in
  \eqref{eqn:use-probability-asymptotics}. Lastly, we recognised that the sum
  over the spousal phenotype values was independent of the genotype laws and
  thus this term is identically one. Revisiting Equation
  \eqref{eqn:mendelian-geno-law}, we find
  \begin{align}
    \mu_s^{(1)}(\bm g^s)
    &= \sum_{\bm g^s_\sigma, \bm g^s_\varphi \in \Omega_m^s}
    \big(1 + O(m^{-1/2})\big)
    \psi(\bm g^s \mid \bm g^s_\sigma, \bm g^s_\varphi)
    \mu_s^{(0)}(\bm g^s_\sigma) \mu_s^{(0)}(\bm g^s_\varphi) \\
    &= \big(1 + O(m^{-1/2})\big)\!
    \expectation_{\mu_s^{(0)} \otimes \mu_s^{(0)}}
    \!\big[\psi(\bm g^s \mid \bm G^s_\sigma, \bm G^s_\varphi)\big] \\
    &= \big(1 + O(m^{-1/2})\big) \mu_s^{(0)}(\bm g^s).
  \end{align}
  In the last step, we identified \(\mu_s^{(0)} = \expectation_{\mu_s^{(0)}
  \otimes \mu_s^{(0)}}[\psi(\bm g^s \mid \bm G^s_\sigma, \bm G^s_\varphi)]\). To
  see why, observe that since the founder population is derived from panmixis,
  their parents have the same panmictic genotype law as their offspring, and
  thus the decomposition of Equation \eqref{eqn:mendelian-geno-law}
  gives rise to
  this identity. With the base case verified, we consider an arbitrary \(t \in
  \mathbf{N}\). Using the same decomposition,
  \begin{align}
    \mu^{(t)}_s(\bm g^s)
    &= \sum_{\bm g^s_\sigma, \bm g^s_\varphi \in \Omega_m^s}
    \big(1 + O(m^{-1/2})\big) \psi(\bm g^s \mid \bm g^s_\sigma,
    \bm g^s_\varphi) \mu_s^{(t - 1)}(\bm g^s_\sigma) \mu_s^{(t - 1)}
    (\bm g^s_\varphi) \\
    &= \sum_{\bm g^s_\sigma, \bm g^s_\varphi \in \Omega_m^s}
    \big(1 + O(m^{-1/2})\big)^{2t - 1} \psi(\bm g^s \mid
    \bm g^s_\sigma, \bm g^s_\varphi) \mu_s^{(0)}(\bm g^s_\sigma)
    \mu_s^{(0)}(\bm g^s_\varphi) \\
    &= \big(1 + O(m^{-1/2})\big)^{2t - 1} \expectation_{\mu_s^{(0)}
    \otimes \mu_s^{(0)}} \big[\psi(\bm g^s \mid \bm G^s_\sigma, \bm
    G^s_\varphi) \big] \label{eqn:use-time-asymptotic-step} \\
    &= \big(1 + O(m^{-1/2})\big) \mu_s^{(0)}(\bm g^s),
  \end{align}
  where dropping the exponent of the asymptotic term of
  Equation~\eqref{eqn:use-time-asymptotic-step} can be justified by expanding
  the term as a binomial summation, where \(t = o(m^{1/2})\) is dominated by the
  \(O(m^{-1/2})\) scaling law.
\end{proof}

The next corollary implies that atomic estimators are asymptotically invariant
to assortment in the sense that their moments are unaffected in the limit \(m
\to \infty\).

\begin{corollary}
  Let \(h: \Omega_m^s \to \mathbf R^{d}\) be a measurable function
  for some local
  segment \(\mathcal G_s \subset \mathcal G\) with \(|\mathcal G_s| = \mathcal
  O(1)\). Then for all \(k \in \mathbf{N}\),
  \begin{equation}
    \expectation_{\mu^{(t)}_s}\!\big[h(\bm G^s)^k\big] =
    \big(1 + O(m^{-1/2})\big)
    \expectation_{\mu^{(0)}_s}\!\big[h(\bm G^s)^k\big].
  \end{equation}
\end{corollary}

\begin{proof}
  This is a direct result of \Cref{thm:assortment-invariance}:
  \begin{equation}
    \expectation_{\mu^{(t)}_s}\!\big[h(\bm G^s)^k\big] =
    \sum_{\bm g^s \in \Omega_m^s} h(\bm g^s) \mu_s^{(t)}(\bm g^s)
    = \sum_{\bm g^s \in \Omega_m^s }\big(1 + O(m^{-1/2})\big)
    h(\bm g^s) \mu_s^{(0)}(\bm g^s)
    = \big(1 + O(m^{-1/2})\big)
    \expectation_{\mu^{(0)}_s}\!\big[h(\bm G^s)^k\big]. \qedhere
  \end{equation}
\end{proof}

\subsection{Convergence to equilibrium and limiting forms}
\label{app:convergence}

\begin{lemma}
  If \(\overline\rho_{\sigma \times \varphi} \theta\|\Sigma_G^{(0)}\|
  \leq 1/4\),
  then there exists a constant \(\mu > 0\) such that iterates of the
  system of Equation \eqref{eqn:multivariate-bulmer} are uniformly bounded above
  by \(\mu I\).
\end{lemma}

\begin{proof}
  We prove \(\|\Sigma\| \leq \mu\) implies \(\|\mathcal B(\Sigma)\| \leq \mu\).
  Notice that we have
  \begin{equation}
    \|\mathcal B(\Sigma)\|
    \leq \frac{1}{2} \bigg(\|\Sigma_G^{(0)}\| + \|\Sigma\| +
    \|\symmatecor\| \|\Lambda\|^2 \|\Sigma_E^{-1/2}\|^2 \|\Sigma\|^2\bigg),
  \end{equation}
  where we used submultiplicativity and the triangle inequality for the norm. If
  \(\|\Sigma\| \leq \mu\), then what we desire is
  \begin{equation}
    \frac{1}{2} \bigg(\|\Sigma_G^{(0)}\| + \mu + \|\symmatecor\| \|\Lambda\|^2
    \|\Sigma_E^{-1/2}\|^2 \mu^2\bigg) \leq \mu.
  \end{equation}
  which gives the quadratic inequality
  \begin{equation}
    \|\Sigma_G^{(0)}\| - \mu + \theta\asmstrength \mu^2 \leq 0.
  \end{equation}
  Thus, the set of admissible radii is the interval
  \begin{equation}
    \mu \in
    \left[
      \frac{1 - \sqrt{1 - 4\theta\asmstrength
      \|\Sigma_G^{(0)}\|}}{2\theta\asmstrength
      \|\Sigma_G^{(0)}\|},
      \frac{1 + \sqrt{1 -
      4\theta\asmstrength\|\Sigma_G^{(0)}\|}}{2\theta\asmstrength
      \|\Sigma_G^{(0)}\|}
    \right].
  \end{equation}
  We choose the lower endpoint of the interval as the smallest radius of the
  invariant set to obtain the desired conclusion.
\end{proof}

\begin{lemma}\label{lemma:env-herit-ratio-is-equivalent}
  For \(j, k \in \mathbf{N} \cup \{\infty\}\),
  \((I - H^{(j)}) \conmult \Sigma_X^{(j)} =
  (I - H^{(k)}) \conmult \Sigma_X^{(k)}\).
\end{lemma}

\begin{proof}
  The result follows from the assumed constancy of \(\Sigma_E\) and
  \begin{align}
    (I - H^{(t)}) \conmult \Sigma_X^{(t)}
    &= I \conmult \Sigma_X^{(t)} - H^{(t)} \conmult \Sigma_X^{(t)} \\
    &= \Sigma_X^{(t)} - \Sigma_\Gamma^{(t)} \\
    &= \Sigma_E. \qedhere
  \end{align}
\end{proof}

\begin{lemma}\label{lemma:form-of-theta}
  Under commutativity of \((\Sigma_\Gamma^{(0)}, \Sigma_E, \symmatecor)\),
  \begin{equation}
    \Theta^{(0)}_{(\infty)}
    = \Sigma_\Gamma^{(0)} \condiv \Sigma_X^{(\infty)}
    = H^{(0)} (I - H^{(0)})^{-1} (I - H^{(\infty)}).
    \label{eqn:total-var-herit-iden}
  \end{equation}
\end{lemma}

\begin{proof}
  Starting from \Cref{lemma:env-herit-ratio-is-equivalent}, use commutativity of
  the components involved to write
  \begin{equation}
    (I - H^{(\infty)}) \Sigma_X^{(\infty)} = (I - H^{(0)}) \Sigma_X^{(0)},
  \end{equation}
  which we solve to find
  \begin{equation}
    \Sigma_X^{(\infty)} = (I - H^{(\infty)})^{-1} (I - H^{(0)}) \Sigma_X^{(0)}.
  \end{equation}
  Now we are able to write
  \begin{align}
    \Theta^{(0)}_{(\infty)}
    &= \Sigma_\Gamma^{(0)} [\Sigma_X^{(\infty)}]^{-1} \\
    &= \Sigma_\Gamma^{(0)} [\Sigma_X^{(0)}]^{-1} (I - H^{(0)})^{-1}
    (I - H^{(\infty)}) \\
    &= H^{(0)} (I - H^{(0)})^{-1} (I - H^{(\infty)}),
  \end{align}
  which is what we wanted to show. \qedhere
\end{proof}

If the population exhibits either entirely or partially negative assortment, in
the sense that any number of eigenvalues of the assorting matrix \(\symmatecor\)
are negative, taking its principal square root requires some care since this
would involve an attempt to take the root of a negative number. Before
proceeding on to the limiting heritability proof, we introduce some utilities to
clarify these cases. For the Hermitian matrix \(A\) with spectrum \(D\), we
write \(|A| = U|D|U^\top\) to denote its counterpart with non-negative spectrum,
intuitively one whose eigenvalues are the absolute values of those of \(D\).
The sign information is preserved using \(\sgn(A) = U \sgn(D) U^\top\) where the
sign function is applied elementwise to the entries of the diagonal matrix
\(D\). Clearly, \(\sgn(A)^2 = I\). Lastly, observe that we may write
\begin{align}
  A
  &= UDU^\top \\
  &= U \sgn(D) |D|U^\top \\
  &= U|D|^{1/2} \sgn(D) |D|^{1/2} U^\top \\
  &= (U|D|^{1/2} U^\top)[U \sgn(D) U^\top](U |D|^{1/2} U^\top) \\
  &= |A|^{1/2} \sgn(A) |A|^{1/2}
\end{align}
or \(\sgn(A) \conmult |A|\) in the notation adopted here.

\begin{theorem}
  If \(\symmatecor\) is either positive or negative definite, or indefinite and
  \((\Sigma_\Gamma^{(0)}, \Sigma_E, \symmatecor)\) commute pairwise, then the
  matrix equation
  \begin{equation}
    H^{(\infty)} \symmatecor H^{(\infty)} - H^{(\infty)} +
    \Theta^{(0)}_{(\infty)} = 0 \label{eqn:dare-equation}
  \end{equation}
  has solution
  \begin{equation}
    H^{(\infty)} = \frac{1}{2} \sgn(\symmatecor)
    \bigg[I - \bigg(I - 4 \sgn(\symmatecor) (\Theta^{(0)}_{(\infty)} \conmult
    |\symmatecor|) \bigg)^{1/2}\bigg] \condiv |\symmatecor|.
  \end{equation}
  If \((\Sigma_\Gamma^{(0)}, \Sigma_E, \symmatecor)\) commute pairwise, then
  \begin{equation}
    H^{(\infty)} = \frac{1}{2} \sgn(\symmatecor)
    \bigg[
      I - \bigg(I - 4 \sgn(\symmatecor)
      |\symmatecor| H^{(0)}(I - H^{(0)})\bigg)^{1/2}
    \bigg]
    \big[|\symmatecor| (I - H^{(0)})\big]^{-1}, \label{eqn:gimelfarb-matrix}
  \end{equation}
  in which case the eigenvalues \(h_{(\infty), j}^2\) of \(H^{(\infty)}\) may be
  identified as
  \begin{equation}
    h_{(\infty), j}^2 =
    \frac{1 - \sqrt{1 - 4 \bar{r}_j h_{(0), j}^2 (1 - h_{(0), j}^2)}}
    {2\bar{r}_j (1 - h_{(0), j}^2)}.
  \end{equation}
\end{theorem}

\begin{proof}
  We begin by considering the general case. First, complete the
  square on the left hand side, giving
  \begin{equation}
    (H^{(\infty)} - \tfrac{1}{2} \symmatecor^{-1}) \symmatecor (H^{(\infty)} -
    \tfrac{1}{2} \symmatecor^{-1})
    = \frac{1}{4} \symmatecor^{-1} - \Theta^{(0)}_{(\infty)},
  \end{equation}
  which we may transcribe into the form
  \begin{equation}
    (H^{(\infty)} - \tfrac{1}{2} \symmatecor^{-1})
    |\symmatecor|^{1/2} \sgn(\symmatecor) |\symmatecor|^{1/2}
    (H^{(\infty)} - \tfrac{1}{2} \symmatecor^{-1})
    = \frac{1}{4} \big(\symmatecor^{-1} - 4\Theta^{(0)}_{(\infty)}\big).
  \end{equation}
  To obtain a quadratic on the left hand side, this suggests left multiplying by
  \(\sgn(\symmatecor) |\symmatecor|^{1/2}\) and right multiplying by
  \(|\symmatecor|^{1/2}\), which gives
  \begin{equation}
    \left[
      \sgn(\symmatecor)
      \bigg(H^{(\infty)} - \tfrac{1}{2} \symmatecor^{-1}\bigg) \conmult
      |\symmatecor|
    \right]^2 = \frac{1}{4}\big(
      I - 4\sgn(\symmatecor) (\Theta^{(0)}_{(\infty)} \conmult |\symmatecor|)
    \big).
  \end{equation}
  Taking the square root of both sides,
  \begin{equation}
    \sgn(\symmatecor)
    \bigg[\big(H^{(\infty)} - \tfrac{1}{2} \symmatecor^{-1}\big) \conmult
    |\symmatecor|\bigg] =
    -\frac{1}{2} \big(I - 4 \sgn(\symmatecor) (\Theta^{(0)}_{(\infty)} \conmult
    |\symmatecor|)\big)^{1/2},
  \end{equation}
  with the negative branch chosen to ensure conformity with the univariate form.
  Expanding the bracketed term on the left hand side and rearranging,
  \begin{equation}
    \sgn(\symmatecor) H^{(\infty)} \conmult |\symmatecor|
    = \frac{1}{2} \bigg[
      I - \big(I - 4 \sgn(\symmatecor) (\Theta^{(0)}_{(\infty)} \conmult
      |\symmatecor|)\big)^{1/2}
    \bigg].
  \end{equation}
  A final left multiplication by \(\sgn(\symmatecor) |\symmatecor|^{-1/2}\) and
  right multiplication by \(|\symmatecor|^{-1/2}\) gives
  \begin{equation}
    H^{(\infty)} = \frac{1}{2} \sgn(\symmatecor) \left[
      I - \big(I - 4 \sgn(\symmatecor) (\Theta^{(0)}_{(\infty)} \conmult
    |\symmatecor|)\big)^{1/2} \right] \condiv |\symmatecor|.
    \label{eqn:heritability-limit-app}
  \end{equation}
  Since positivity of the discriminant \(I - 4 \sgn(\symmatecor)(
  \Theta^{(0)}_{(\infty)} \conmult |\symmatecor|)\) is required for the matrix
  square root to exist, we derive the condition
  \(\sgn(\symmatecor)(\Theta^{(0)}_{(\infty)} \conmult |\symmatecor|) \preceq
  \frac{1}{4} I\). In the case where \(\symmatecor\) is indefinite,
  \(\sgn(\symmatecor) (\Theta_{(\infty)}^{(0)} \conmult |\symmatecor|)\) may
  not be symmetric, in which case the above form is ill-posed. This is only
  problematic in the indefinite regime, since otherwise \(\sgn(\symmatecor) =
  \pm I\). Under pairwise commutativity of \((\Sigma_\Gamma^{(0)}, \Sigma_E,
  \symmatecor)\), however, all of the factors of \(\sgn(\symmatecor)
  (\Theta_{(\infty)}^{(0)} \conmult |\symmatecor|)\) commute, which is
  equivalent to the matrix being symmetric. Thus, a closed form is available to
  indefinite assorting matrices only in the commuting case. Turning now to this
  special case more generally, we revisit Equation \eqref{eqn:dare-equation},
  which by \Cref{lemma:form-of-theta} may be written
  \begin{equation}
    \symmatecor [H^{(\infty)}]^2 - H^{(\infty)} + H^{(0)} (I -
    H^{(0)})^{-1} (I -
    H^{(\infty)}) = 0.
  \end{equation}
  Multiplying by \(I - H^{(0)}\) and simplifying, this is
  \begin{equation}
    \symmatecor (I - H^{(0)}) [H^{(\infty)}]^2 - H^{(\infty)} + H^{(0)} = 0.
  \end{equation}
  Solving this equation may be done in roughly the same spirit as that of
  solving a scalar quadratic owing to the commutativity of the matrices
  involved. Pursuing such an endeavour gives the desired form of Equation
  \eqref{eqn:gimelfarb-matrix}. Lastly, we may perform the
  eigendecomposition of
  all the matrices in said form with respect to the shared eigenbasis \(U\),
  writing \(\mathcal{R}\), \(\mathcal H^{(0)}\), and \(\mathcal H^{(\infty)}\)
  to denote the diagonal spectra of \(\symmatecor\), \(H^{(0)}\), and
  \(H^{(\infty)}\), respectively. Then, Equation \eqref{eqn:gimelfarb-matrix}
  takes the form
  \begin{equation}
    \mathcal H^{(\infty)} = \frac{1}{2} U \bigg[I - \bigg(I - 4\mathcal{R}
        \mathcal{H}^{(0)}(I - \mathcal
    H^{(0)})\bigg)^{1/2}\bigg]\big[\mathcal{R}(I -
    \mathcal{H}^{(0)})\big]^{-1}U^\top.
  \end{equation}
  Reading off the diagonals of the above gives the desired identity.
\end{proof}

\end{document}